\documentclass[EJP]{ejpecp}
\pdfoutput=1

\AUTHORS{%
  Leandro~P.~R.~Pimentel\footnote{Federal University of Rio de Janeiro, Brazil.
    \EMAIL{leandro@im.ufrj.br}}
  \and 
  Rapha\"el~Rossignol\footnote{Universit\'e Joseph Fourier Grenoble 1, France. \EMAIL{raphael.rossignol@ujf-grenoble.fr}} }

\SHORTTITLE{Greedy Polyominoes and first-passage times}

\TITLE{Greedy Polyominoes
  and first-passage times on random Voronoi tilings\thanks{R. Rossignol was partially suported by the Swiss National Science
   Foundation grants 200021-1036251/1 and 200020-112316/1.}\thanks{L. Pimentel was partially supported by grants from the Swiss National Science Foundation, Funda\c{c}\~ao de Amparo a Pesquisa do Estado de S\~{a}o Paulo and The Netherlands Organisation for Scientific Research.}} 


\AMSSUBJ{60K35} 
\AMSSUBJSECONDARY{60D05} 

\SUBMITTED{January 12, 2008} 
\ACCEPTED{January 31, 2012} 

\ARXIVID{0811.0308}  

\VOLUME{0}
\YEAR{2012}
\PAPERNUM{0}
\DOI{vVOL-PID}

\ABSTRACT{
Let $\N$ be distributed as a Poisson random set on $\RR^d$, $d\geq 2$, with
intensity comparable to the Lebesgue measure. Consider the
Voronoi tiling of $\RR^d$, $\{C_v\}_{v\in \N}$, where $C_v$ is
composed of points $\x\in\RR^d$ that are closer to $v\in\N$ than to
any other $v'\in\N$.  A polyomino $\P$ of size $n$ is a connected
union (in the usual $\RR^d$ topological sense) of $n$ tiles, and we denote by $\Pi_n$ the collection of all
polyominos $\P$ of size $n$ containing the origin. Assume that the weight of a Voronoi tile $C_v$ is given by $F(C_v)$, where $F$ is a nonnegative functional on Voronoi tiles. In this paper we investigate for some functionals $F$, mainly when $F(C_v)$ is a polynomial function of the number of faces of $C_v$,  the tail behavior of the maximal weight among polyominoes in
$\Pi_n$: $F_n=F_n(\N):=\max_{\P\in\Pi_n} \sum_{v\in \P}
F(C_v)$. Next we apply our results to study self-avoiding paths, first-passage percolation models and the stabbing number on the dual graph, named the Delaunay triangulation. As the main application we show that first passage percolation has at most linear variance.
}

\newcommand{\0}{\mathbf{0}}
\newcommand{\x}{\mathbf{x}}
\renewcommand{\v}{\mathbf{v}}
\newcommand{\z}{\mathbf{z}}
\newcommand{\N}{\mathcal{N}}
\newcommand{\D}{\mathcal{D}}
\renewcommand{\P}{\mathcal{P}}
\newcommand{\B}{\mathbf{B}}

\newcommand{\A}{\mathbf{A}}
\newcommand{\SA}{\mathcal{SA}}

\newcommand{\V}{\mathbf{V}}
\newcommand{\G}{\mathcal{G}}
\newcommand{\E}{\mathcal{E}}

\newcommand{\calC}{\mathcal{C}}
\newcommand{\Cl}{\mathbf{Cl}}

\newcommand{\Adh}{Ad}

\renewcommand{\Im}{{Bl}}

\newcommand{\C}{\mathbf{C}}
\newcommand{\PP}{\mathbb{P}}
\newcommand{\GG}{\mathbb{G}}
\newcommand{\RR}{\mathbb{R}}

\newcommand{\EE}{\mathbb{E}}
\newcommand{\ER}{\mathcal{E}_{d}} 
\newcommand{\NN}{\mathbb{N}}

\newcommand{\ZZ}{\mathbb{Z}}

\newcommand{\Var}{\mathsf{Var}}
\newcommand{\II}{\mbox{ 1\hskip -.29em I}}
\newcommand{\sfrac}[2]{\kern.1em
        \raise.5ex\hbox{$#1$}\kern-.1em
        /\kern-.15em\lower.25ex\hbox{$#2$}}

\newcommand{\bdes}{\begin{description}}
\newcommand{\edes}{\end{description}}

\newcommand{\eps}{\varepsilon}
\newcommand{\bex}{\begin{example} \em }
\newcommand{\eex}{\end{example} }
\newcommand{\brem}{\begin{remark} \em}
\newcommand{\erem}{\end{remark} }

\newcommand{\inter}[1]{{\overset{o}{#1}}}

\begin{document}



\section{Introduction}
Greedy animals on $\ZZ^d$ have been studied notably in
\cite{Coxetal93}. Imagine that positive weights, or awards, are placed
on all vertices of $\ZZ^d$. A greedy animal of size $n$ is a connected
subset of $n$ vertices, containing the origin, and which catches the
maximum amount of awards. When these awards are random, i.i.d, it is
shown in \cite{Coxetal93} that the total award collected by a greedy
animal grows at most linearly in $n$ if the tail of the award is not too thick. This is shown in a rather strong sense, giving deviations
inequalities decaying at a rate $n^{-\log(n)^\alpha}$, for some $\alpha$ (see Proposition 1 in \cite{Coxetal93}). This result has already been proved useful for the study
of percolation and (dependent) First Passage Percolation, see
\cite{FontesNewman93}.

\begin{figure}[htb]
\begin{center}
\includegraphics[width=0.5\textwidth]{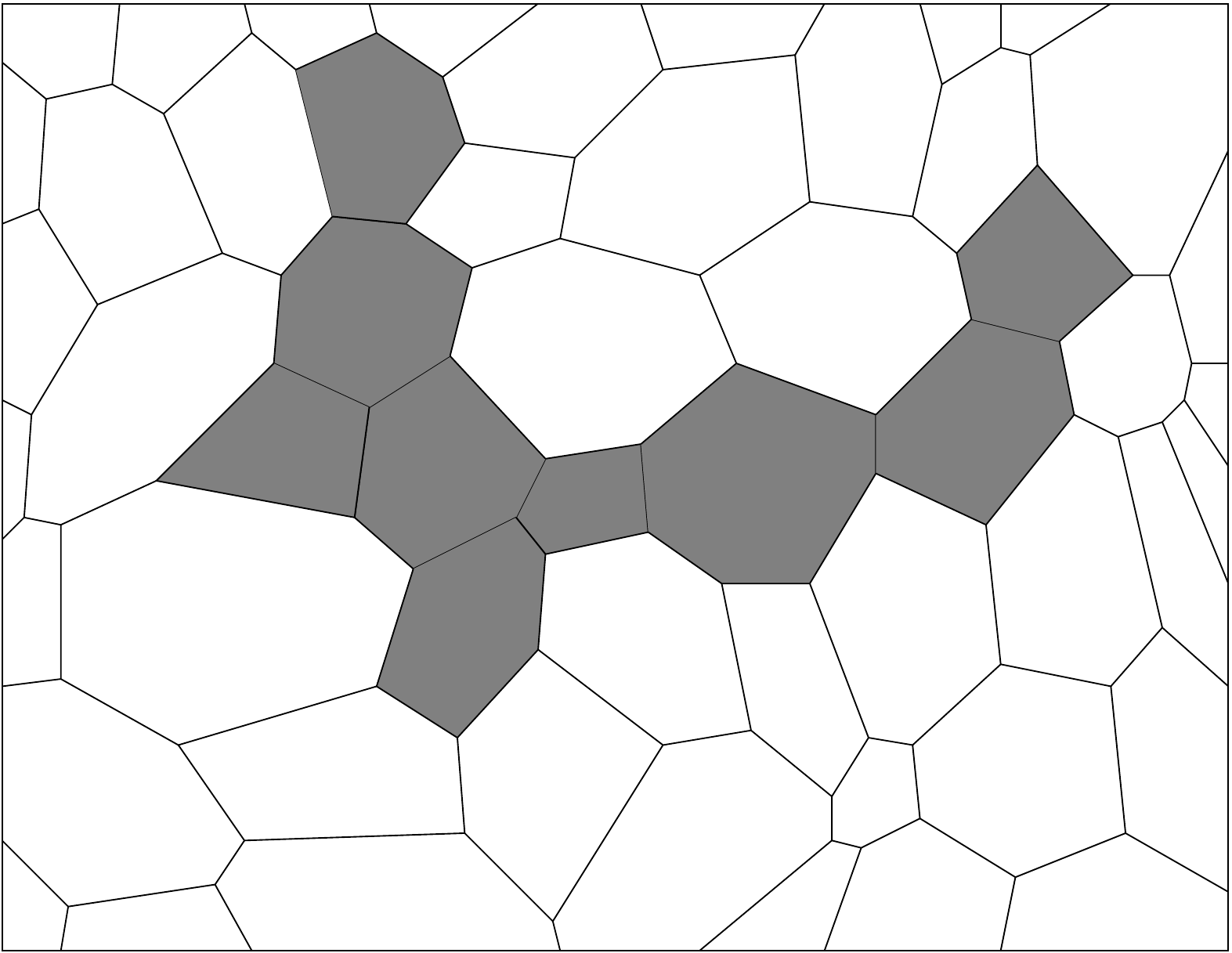}
\caption{A two-dimensional Voronoi polyomino of size $n=9$.}
\label{poly}
\end{center}
\end{figure}

The aim of the present paper is to extend the results of \cite{Coxetal93} to
some \emph{ greedy polyominoes on the Poisson-Voronoi tiling}. The precise definitions will be given in section \ref{sec:definitions}, but let us explain what we mean. First, in this paper, the Poisson-Voronoi tiling is the Voronoi tiling based on a Poisson random set on $\RR^d$ with intensity comparable to the Lebesgue measure. A polyomino of size $n$ on this tiling is a connected union of $n$ Voronoi tiles, cf. Figure \ref{poly}. Then, we let $f$ be some fixed function from $\NN$ to $\RR^+$ (a \emph{weight function}), and put down on each Voronoi tile an award depending on the number of its faces: a Voronoi tile with $r$ faces receives an award equal to $f(r)$. Finally, a polyomino of size $n$ is greedy if it contains the origin and catches the maximum amount of awards.

Our main result, Theorem \ref{majotheo} below, gives a deviation inequality for the total award collected by a greedy polyomino on the Poisson-Voronoi tiling. Of course, the rate of decay depends on the precise weight function $f$ defining the award. To give an idea of our results, when $f(r)=r^k$, $k\geq1$, we obtain the following deviation inequality.
\begin{corollary}
\label{coromajotheo}
Let $k\geq 0$ and denote by $F_n(k)$ the total award collected by a greedy polyomino on the Poisson-Voronoi tiling, with award on a tile being equal to the $k$-th power of the number of faces of the tile. Then, there are constants $z_0$ and $K$ such that for any $n\geq 1$, and any $z\geq z_0$,
$$\PP\left(\frac{1}{n}F_n(k)\geq z\right)\leq e^{-K(nz)^\frac{1}{k+2}}\;.$$
\end{corollary}

Our results may be useful to control the geometry of the Poisson-Voronoi tiling, and this may be best viewed through the facial dual of the Poisson-Voronoi tiling, which is called the \emph{Poisson-Delaunay graph}. Notice that through this duality, a greedy polyomino on the Poisson-Voronoi tiling with weight function $f$ is simply a greedy animal on the Delaunay graph where the award on a vertex of degree $k$ equals $f(k)$. We shall give three applications of our results on greedy polyominoes to some geometric problems on the \emph{Poisson-Delaunay graph}. 

The first application is an estimate on the tail probability of the number of
self-avoiding paths of length $n$ starting from the origin. 

The second application concerns First Passage Percolation on the Poisson-Delaunay graph.  We prove that the first
passage time in First Passage Percolation has (at most) linear variance on the
Poisson-Delaunay graph. 

Finally, our third application concerns tail estimates on the stabbing number, which is defined in \cite{AddarioberrySarkar05} as the maximum number of Delaunay cells that intersect a single line in the cube $[0,n]^d$. This is an important quantity, as can be seen in \cite{AddarioberrySarkar05} where it is the crucial ingredient to derive the transience or recurence properties of the simple random walk on the Poisson-Delaunay graph. We shall obtain an exponential deviation inequality for the stabbing number of the Poisson-Delaunay graph.

The rest of the paper is organized as follows. In section \ref{sec:definitions}, we give the precise definitions and state our main results concerning greedy polynominoes on the Poisson-Voronoi tiling (and greedy animals on the Poisson-Delaunay graph). Section \ref{sec:proofs} is devoted to the proofs of these main results, which is  based on a renormalization argument from percolation theory and an adaptation of the chaining technique of \cite{Coxetal93}. This last step being rather technical, its proof is postponed to the appendix. Finally, the three applications on the Poisson-Delaunay graph are the matter of section \ref{sec:applications}.

\section{Definitions and main results}
\label{sec:definitions}
In this section, we give all the needed formal definitions and state precisely our main results on greedy polyominoes, Theorem \ref{majotheo} and Proposition \ref{prop:Gamma}.

\subsection{The Poisson-Voronoi tiling and polyominoes.}
In the whole paper, we suppose that $d\geq 2$. To any locally finite subset $\N$ of $\RR^d$ one can associate a  collection of subsets of $\RR^d$ whose union is $\RR^d$. To each point $v\in \N$
corresponds a polygonal region $C_v$, the \emph{Voronoi tile} (or
cell) at $v$, consisting of the set of points of $\RR^d$ which are
closer to $v$ than to any other $v'\in \N$. Closer is understood here in the large sense, and
this collection is not a partition, but the set of points which belong to
more than one Voronoi tile has Lebesgue measure $0$. The collection
$\{C_v\}_{v\in \N}$ is called the \emph{Voronoi tiling} (or
tessellation) of the plane based on $\N$. From now on, $\N$ is
understood to be distributed like a Poisson random set on $\RR^d$ with
intensity measure $\mu$. We shall always suppose that $\mu$ is
comparable to Lebesgue's measure on $\RR^d$, $\lambda_d$, in the sense
that there exists a positive constant $c_\mu$ such that for every
Lebesgue-measurable subset $A$ of $\RR^d$:
\begin{equation}
\label{eq:comparaison}
\frac{1}{c_\mu}\lambda_d(A)\leq\mu(A)\leq c_\mu\lambda_d(A)\;.
\end{equation}
For each positive integer number $n\geq 1$, a \emph{Voronoi
polyomino} $\P$ of size $n$ is a connected union of $n$
Voronoi tiles (Figure \ref{poly}). Notice that with probability one,
when two Voronoi tiles are connected, they share a $(d-1)$-dimensional face.  We denote by $\Pi_n$ the set of
all polyominoes $\P$ of size $n$ and such that the origin $0$ belongs to
 $\P$.

Assume that the ``weight'' of a Voronoi tile $C_v$ is given by
$f(d_{\N}(v))$, where $d_{\N}(v)$ is the number of $(d-1)$-dimensional faces of $C_v$ and
$f$ is a nondecreasing function from $[1,\infty)$ to
$[1,\infty)$. In this way we define a random weight functional on polyominoes by
$$
F(f,\N,\P):=\sum_{v\in \P} f(d_{\N}(v))\,.
$$
The maximal weight among polyominoes in $\Pi_n$ is:
\begin{equation}
\label{eq:greedy}
F_n(f,\N):=\max_{\P\in\Pi_n} F(f,\N,\P)\;.
\end{equation}
A greedy Voronoi polyomino $\P_n$ is a Voronoi polyomino that attains
the maximum in the definition of $F_n$: $F(f,\N,\P_n)=F_n(f,\N)$.

To state our main theorem, we require the following notation:
\begin{definition}\label{functions} Let $g^{-1}$ denote the pseudo-inverse of any strictly increasing function $g$ from $[1,\infty)$ to
$[1,\infty)$:
$$g^{-1}(u)=\sup\{u'\in [1,\infty)\mbox{ s.t. }g(u')< u\}\;.$$
A \emph{ weight function} $f$ is a nondecreasing function from $[1,\infty)$ to
$[1,\infty)$. Given a weight function $f$, we define two functions $\hat f$ and $f^*$ as follows:
$$\hat f=(x\mapsto xf(x))^{-1}\;,$$
and
$$f^*=\hat f\circ (y\mapsto y\hat f(y))^{-1}\;.$$
In addition, for each $c\in(0,\infty)$, we say that $f$ is $c$-nice if
$$
\liminf_{y\to\infty}\frac{\hat f(y)}{\log y}\geq c\;.
$$
\end{definition}
Examples of $c$-nice weight functions $f$, and various $\hat f$, $f^*$ are given after Theorem \ref{majotheo}.
\begin{theorem}
\label{majotheo} Let $f$ be a non-decreasing weight function. There exists constants $z_1,c_1,c_2$ and $c_3$ in $(0,\infty)$ such that if $f$ is $c_1$-nice then for all $z\geq z_1$ and for all
$n\geq 1$:
$$
\PP\big(F_n(f,\N)>nz\big)\leq \exp\Big\{-c_2 f^*(c_3nz)\Big\}\;.
$$
In particular, there also exists a constant $c_4\in(0,\infty)$ such
that
$$
0\leq f(3)\leq\liminf_{n\to\infty}\EE\Big(
\frac{F_n(f,\N)}{n}\Big)\leq\limsup_{n\to\infty}\EE\Big(
\frac{F_n(f,\N)}{n}\Big)<c_4\;.
$$
\end{theorem}
\begin{example}
If $f(x)=x^k$ with $k\geq 0$, then $\hat f(u)=u^{\frac{1}{k+1}}$ and $f^*(x)=x^{\frac{1}{k+2}}$. This weight function $f$ is $c$-nice for any $c>0$, and this gives Corollary \ref{coromajotheo}.
\end{example}
\begin{example}
Let $f(x)=\exp\left\{(x/K)^\alpha\right\}$, with $K>0$ and $\alpha>0$. Then as $x$ goes to infinity, $\hat f(u)\sim K(\log u)^{1/\alpha}$ and $f^*(x)\sim K(\log x)^{1/\alpha}$. Thus if $\alpha<1$, $f$  is $c$-nice for any $c>0$, and Theorem \ref{majotheo} gives 
$$\PP\big(F_n(f,\N)>nz\big)\leq \exp\Big\{-C(\log(nz))^\frac{1}{\alpha}\Big\}\,,$$
for some constant $C>0$ and $z$ large enough. Also, when $\alpha=1$, $f$ is $c$-nice for any $c\in (0, K]$. Thus, if  $K$ is large enough (at least the constant $c_1$ of Theorem \ref{majotheo}), then:
$$\PP\big(F_n(f,\N)>nz\big)\leq \exp\Big\{-C\log(nz)\Big\}\,,$$
for some constant $C>0$ and $z$ large enough. 
\end{example}

\subsection{The Delaunay graph.} 
\label{subsec:defDelaunay}
\begin{figure}[htb]
\begin{center}
\includegraphics[width=0.5\textwidth]{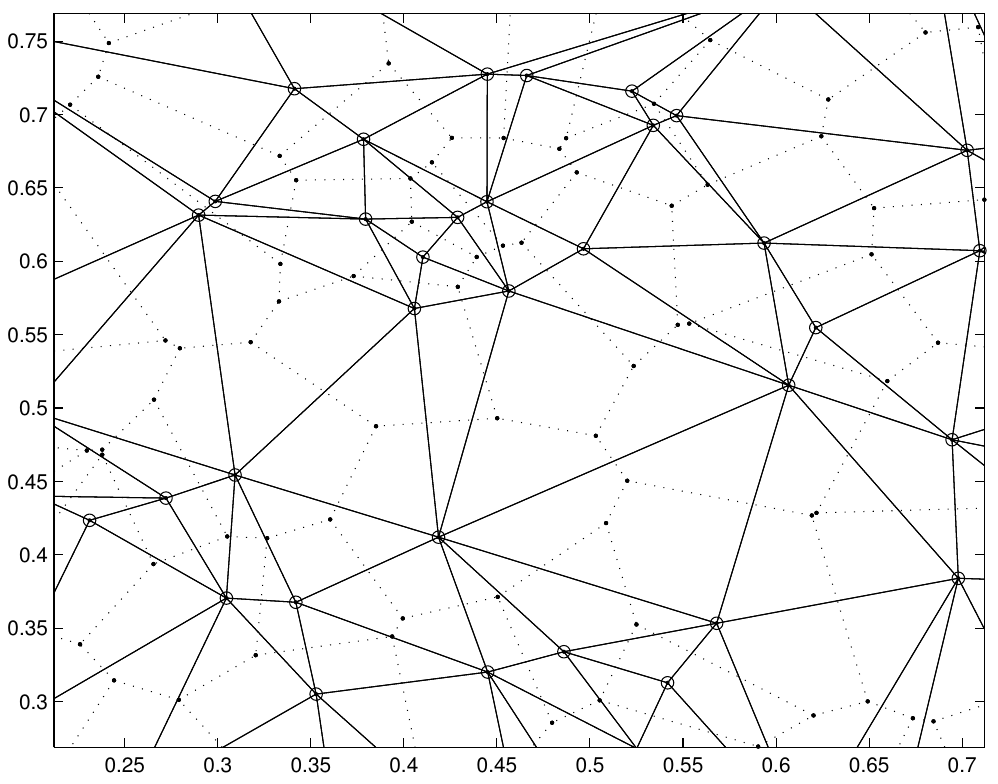}
\caption{The Voronoi tiling (dashed lines) and the Delaunay triangulation (solid lines) in dimension $d=2$.}
\label{vor}
\end{center}
\end{figure}

An important graph for the study of a Voronoi tiling is its facial dual, the \emph{Delaunay graph} based on $\N$. This graph, denoted by $\D(\N)$ is an unoriented graph embedded in $\RR^d$ which has vertex set $\N$ and edges $\{u,v\}$ every time $C_u$ and $C_v$ share a $(d-1)$-dimensional face (Figure \ref{vor}). We remark that, for our Poisson random set, almost surely no $d+1$ points are on the same hyperplane and no $d+2$ points are on the same hypersphere, which makes the Delaunay graph a well defined triangulation. This triangulation divides $\RR^d$ into bounded simplices called  \emph{Delaunay cells}. For a Delaunay cell $\Delta$, let $B(\Delta)$ denote the closed circumball of $\Delta$. An important property which we shall use several times is that for each Delaunay cell $\Delta$ no point in $\N$ lies in the interior of $B(\Delta)$. Polyominoes on the Voronoi tiling correspond to connected (in the graph topology) subsets of the Delaunay graph. Also, the number of faces of a Voronoi tile $C_v$, which we denoted by $d_\N(v)$, is simply the degree of $v$ in  $\D(\N)$. 

We shall prove an easier variant of Theorem \ref{majotheo}, which
will be useful in the applications. Define $\Omega$ to be the set of locally finite subsets of $\RR^d$. Then, for any $\omega\in\Omega$ and any subgraph $\phi$ of $\mathcal{D}(\omega)$, we define $\Gamma(\omega,\phi)$ to be the set of
points in $\RR^d$ whose addition to $\omega$ ``perturbs'' $\phi$: 
$$\Gamma(\omega,\phi)=\{x\in\RR^d\mbox{ s.t. }\phi\not\subset
\mathcal{D}(\omega\cup \{x\})\}\;.$$
We get the following result for the maximal size of
$\Gamma(\omega,\phi)$ when $\phi$ belongs to $\SA_n$, the set of
self-avoiding paths on $\D(\N)$ starting from $v(0)$ and of size $n$, where $v(0)$ is the a.s. unique $v\in\N$ s.t. $0\in C_v$: 
\begin{proposition}
\label{prop:Gamma}
There are constants $z_1$ and $C_1$ such that for every $z\geq z_1$,
and for every $n\geq 0$,
$$
\PP\big(\max_{\phi\in\SA_n} \mu(\Gamma(\N,\phi))>nz\big)\leq e^{-C_1nz}\;.
$$
\end{proposition}

\section{Proofs of the main results}
\label{sec:proofs}

There are three main issues when considering the total award of a greedy polyomino as in (\ref{eq:greedy}). The first issue lies in the nature of the weights: they are a function of the number of faces of a tile or equivalently, the degree of a vertex in the Delaunay graph. Notice that it is known how to control this quantity for the ``typical'' cell, in the Palm sense (cf. \cite{Calka03} for instance) but we are not in this ``typical'' case. Of course, we know really well how to count the number of points of the Poisson random set in a fixed area, but we do not know that well how far away we have to expect the neighbours of these points to lie. The idea to answer this problem is to use a renormalization trick from percolation theory. More precisely, we shall consider a box in $\RR^d$ large enough so that it contains with a ``large enough'' probability some configuration of points which prevents a Delaunay cell to cross it
completely. This ``large enough''probability corresponds to a
percolation threshold: we need that the ``bad boxes'' (those who can
be crossed) do not percolate. This will allow us to control the degree of a vertex in a Polyomino, bounding it from above by the number of points inside a cluster of ``bad boxes'' containing the vertex. In section \ref{comp}, we explain this renormalization by showing how to cover animals with boxes and
clusters of boxes. 

A second issue lies in the fact that the degrees of two vertices are dependent. This will be solved inside the renormalization trick, at the price of counting many times the same number of points in a cluster of ``bad boxes''.

The third issue is to handle the fact that the supremum in (\ref{eq:greedy}) is over a potentially high (i.e exponential in $n$) number of polyominoes of size $n$. This will be handled through the chaining technique due to \cite{Coxetal93}. We shall state the corresponding adaptation in section \ref{subsec:chaininglemma}.

With these tools in hand, we shall prove our main results, Theorem \ref{majotheo} in section \ref{subsec:proofmain} and Proposition \ref{prop:Gamma}
in section \ref{subsec:proofgamma}.

\subsection{The renormalization trick: comparison with site percolation and lattice animals}
\label{comp}

We denote by $\GG_d$ the $d$-dimensional lattice with vertex set
$\ZZ^d$ and edge set composed by pairs $(\z,\z')$ such that
$$|\z-\z'|_{\infty}=\max_{j=1,\dots,d}|\z(j)-\z'(j)|= 1\,.$$  
A box $B$ in $\RR^d$ is any set of the form 
$$B=\x+[0,L)^d\,$$
for some $\x\in\RR^d$ and $L\geq 0$. We shall often use boxes centered at points of some lattice embedded in $\RR^d$, so we define, for any $z$ in $\RR^d$
$$B_\z=\z+[-1/2,1/2)^d\,.$$
Now, we define the notions of
nice and good boxes. These are notions depending
on a set $\N$ of points in $\RR^d$, which basically say that there
are enough points of $\N$ inside the box so that a Voronoi tile of a
point of $\N$ outside the box cannot ``cross'' the box. This will be
given a precise sense in Lemma \ref{lemgoodbox}. 

\begin{definition}
\label{defi:goodbox}
Define the following even integer:
$$\alpha_d=2(4\lceil\sqrt{d}\rceil+2)\;.$$
Let $\N$ be a set of points in $\RR^d$.
We say that a box $B\subseteq\RR^d$ is $\N$-nice if, cutting it
regularly into $\alpha_d^d$ sub-boxes, each one of these boxes
contains at least one point of the set $\N$. 
We say that a box $B$ is $\N$-good if, cutting it
regularly into $(3\alpha_d)^d$ sub-boxes, each one of these boxes
contains at least one point of the set $\N$. When $\N$ is
understood, we shall simply say that $B$ is \emph{nice} (resp. \emph{good}) if it is
$\N$-nice (resp. $\N$-good). We say a box is
$\N$-\emph{ugly} (resp. $\N$-\emph{bad}) if it is not $\N$-nice (resp. not $\N$-good).
\end{definition}
Notice that if a box is good, then it is nice, and thus if a box is ugly, then it is bad. To count the number of points of $\N$ which are in a subset $A$ of
$\RR^d$, we define:
$$
|A|_\N=|A\cap \N|\;.
$$

An animal $\A$ in $\GG_d$ is a finite and connected subset of $\GG_d$. To each bounded and connected subset $A$ of
$\RR^d$, we associate an animal in $\GG_d$:
$$
\A(A)=\{\z\in\GG_d\mbox{ s.t.
}A\cap B_\z\not = \emptyset\}\;.
$$
If $\V$ is a set of vertices of $\ZZ_d$, define the border $\partial \V$ of $\V$ as the
set of vertices which are not in $\V$ but have a $\GG^d$-neighbor in
$\V$. We also define the following subsets of $\RR^d$:
$$\Im(\V):=\bigcup_{\z\in\V}B_\z\;.$$
$$\Adh(\V):=\left\{x\in\RR^d\mbox{ s.t. }\inf_{y\in \Im(\V)}\|x-y\|_\infty<\frac{1}{2}\right\}\;.$$

A site percolation scheme on $\GG_d$ is defined by:
$$X_\z=\II_{B_\z\;\mathrm{is}\;\N-bad}\,,\,\,\mbox{ for }\z\in
\ZZ^d\;.
$$
$X=(X_\z)_{\z\in\ZZ_d}$ is a collection of independent Bernoulli random
variables. A vertex $\z$  is said to be \emph{bad} when
$X_\z=1$. A \emph{bad cluster} is then a maximal
 connected subset $\C$ of bad vertices of $\GG_d$. Similarly, we
 define ugly vertices and ugly clusters. For any set of
 vertices $\V$, we define by $\Cl(\V)$ the collection of all the bad
 clusters intersecting $\V$. We shall make a slight abuse of notation
 by writing $\Cl(\z)$ to be the bad cluster containing $\z$, if
 there is any (otherwise, we let $\Cl(\z)=\emptyset$).

We now want to show that the Delaunay cells cannot cross the boundaries of an ugly (or bad) cluster. Since the argument will be used at other places, we state a more general lemma first.

\begin{lemma}
\label{lemtruegoodbox}
Let $\C$ be a non-empty connected set of vertices in $\GG_d$. define $\partial \C$ to be the exterior boundary of $\C$:
$$\partial \C=\{\z\in\GG_d\mbox{ s.t. }z\not\in \C\mbox{ and }\exists x\in\C\mbox{ s.t. }|\z-\z'|_{\infty}=1\}\,.$$
Suppose that $\partial \C$ is composed of nice vertices and let $B'$ be a ball in $\RR^d$ such that:
$$B'\cap \Im(\C)\not=\emptyset\quad\mbox{and}\quad\inter{B'}\cap \N=\emptyset\,.$$
Then,
$$B'\subset \Adh(\C)\;.$$
\end{lemma}

\begin{proof}
The idea of the proof is essentially the same as in Lemma 2.1 from
\cite{Pimentel05}. We proceed by contradiction. Suppose that 
$$B'\cap \Im(\C)\not=\emptyset,\quad \inter{B'}\cap \N=\emptyset\quad\mbox{and}\quad B'\not\subset \Adh(\C)\,.$$
Let us define:
$$\partial^{\infty} \C=\Adh(\C)\setminus\Im(\C)\;.$$
This implies that
one may find $x_1$, $x_2$ in $B'$ and $y$ in $\inter{B'}$, and a box $B$ of side
length $\frac{1}{2}$ such that:
\begin{equation}
\label{eq:cross}
\begin{array}{l}y\in[x_1,x_2],\;y\in \mbox{central}(B),\;x_1\not\in
  B,\;x_2\not\in B\;,\\
\mbox{and  cutting }B\mbox{ regularly into
}(\alpha_d/2)^d\mbox{ sub-boxes},\\
\mbox{each sub-box
contains at least one point of }\N\;,
\end{array}
\end{equation}
where $\mbox{central}(B)$ is the only sub-box of $B$ containing the
center of $B$ when one cuts $B$ regularly into $(\alpha_d/2)^d$
sub-boxes. The ball $B'$ necessarily has diameter greater than $\|x_1-x_2\|_2$, 
which is larger than $1/2$. Thus, there is a ball $B'$ of diameter larger
than $1/2$, which contains $y$ and whose interior does not contain any point of
$\N$. But one may see that any ball of diameter strictly larger
than $2\sqrt{d}/\alpha_d$ necessarily contains at least one of the
sub-boxes of side length $1/\alpha_d$. Notice that any ball of diameter strictly smaller than
$\frac{(\alpha_d-1)}{2\alpha_d}$ which contains a point in $\mbox{central}(B)$ is totally included in $B$. Since 
$$2\sqrt{d}/\alpha_d< \frac{(\alpha_d-1)}{2\alpha_d}\;,$$
and since the sub-boxes of $B$ of side length $1/\alpha_d$ contain a point of $\N$, we deduce that $B'$ contains a sub-box of
$B$ which contains a point of $\N$, whence a contradiction. 
\end{proof}
\begin{lemma}
\label{lemgoodbox}
Assume that $\C$ is an ugly cluster in $\GG_d$. Let $\Delta$ be any Delaunay cell of $\D(\N)$. Then, 
$$\Delta\cap \Im(\C)\not=\emptyset\Rightarrow \Delta\subset
\Adh(\C)\;.$$
The same holds for bad clusters.
\end{lemma}

\begin{proof}
Suppose that
$\Delta\cap \Im(\C)\not =\emptyset$. Notice that the exterior boundary of an ugly cluster is composed of nice vertices. Since $\Delta$ is a Delaunay cell, the circumball of $\Delta$, $B(\Delta)$, is a ball containing $\Delta$ and whose interior does not contain any point of $\N$. Thus, one may apply Lemma \ref{lemtruegoodbox} to deduce that $B(\Delta)$, and hence $\Delta$, is included in $\Adh(\C)$. 

The same proof holds for bad clusters, since the boundary of a bad cluster is composed of good vertices, which are nice
vertices.
\end{proof}

The renormalization trick is essentially to transfer the problem of counting the degrees in
our functional to counting the points in some boxes or clusters of bad boxes which cover the animals. However, to have that bad boxes occur with small probability we need to rescale our initial boxes. For $\z$ in $\ZZ^d$ and $r>0$, we define
$$B_\z^{r}=r\z+[-r/2,r/2)^d\,.$$
In this rescaled setup, good and nice boxes are defined according to Definition \ref{defi:goodbox}. In order to cover random polyominoes properly we need to translate the percolation scheme as follows: 
$$
\forall i=1,\ldots ,3^d,\;X_\z^{r,i}=\II_{r\vec{f}_i/3+B_\z^{r}\;\mathrm{is}\;\N-bad},\;\forall \z\in
\ZZ^d\;,
$$
where $\vec{f}_1=0$ and $\vec{f}_2,\dots,\vec{f}_{3^d}$ are the
neighbors of $\0$ in $\GG_d$. For each $i=1,\dots,3^d$ and $r>0$, $X^{r,i}=(X_\z^{r,i})_{\z\in\ZZ_d}$ is a collection of independent Bernoulli random
variables. Of course, the
comparison of $\mu$ to Lebesgue's measure (\ref{eq:comparaison}) implies that
$\mu(B_\z^{r,i})$ goes to infinity when $r$ goes to infinity. Thus,
$$\lim_{r\to\infty}\sup_i\sup_\z\PP(X_\z^{r,i}=1)=0\,.$$ 
Every time we address to the $X^{r,i}$ percolation scheme we label a name, or a variable, with a $(r,i)$. For instance, a vertex $\z$  is said to be $(r,i)$-\emph{bad} (or
simply \emph{bad} when $r$ and $i$ are implicit) when $X_\z^{r,i}=1$. To each  bounded and connected subset $A$ of
$\RR^d$, we also associate $3^d$ different animals in $\GG_d$:
$$
\forall i=1,\ldots ,3^d,\;\A^{r,i}(A)=\{\z\in\GG_d\mbox{ s.t.
}A\cap (r\vec{f}_i/3+B_\z^r)\not = \emptyset\}\;;
$$
and $3^d$ subsets of $\RR^d$:
$$\forall i=1,\ldots
,3^d,\;\Im^{r,i}(\V):=\bigcup_{\z\in\V}(B_\z^r+r\vec{f}_i/3)\;,$$
and 
$$\forall i=1,\ldots
,3^d,\;\Adh^{r,i}(\V):=\left\{x\in\RR^d\mbox{ s.t. }\inf_{y\in \Im^{r,i}(\V)}\|x-y\|_\infty<\frac{r}{2}\right\}\;.$$
Of course, Lemma \ref{lemgoodbox} still holds for the $X^{r,i}$ setup. Let $\Delta$ be any Delaunay cell of $\D(\N)$. If $\C$ is an $(r,i)$-ugly cluster in $\GG_d$ then 
$$\Delta\cap \Im^{r,i}(\C)\not=\emptyset\Rightarrow \Delta\subset
\Adh^{r,i}(\C)\;.$$

The following lemma makes precise the renormalization trick.
\begin{lemma}
\label{lemfrontiere} Let $\G$ be a finite collection of bounded and connected subsets of $\RR^d$. Then, for any positive real number $r,t>0$.
\begin{eqnarray*}
\PP\left\lbrack\sup_{\gamma\in \G}\sum_{v\in \N\cap\gamma}f(d_\N(v))>
t\right\rbrack&&\leq\\
\sum_{i=1}^{3^d}&&\PP\left\lbrack\sup_{\gamma\in\G}\sum_{\C\in\Cl^{3r,i}(\A^{3r,i}(\gamma))}|\Adh^{3r,i}(\C)|_{\N}f(|\Adh^{3r,i}(\C)|_\N)>\frac{t}{3^d.2} \right\rbrack\\
+\sum_{i=1}^{3^d}&&\PP\left\lbrack\sup_{\gamma\in\G}\sum_{\z\in \A^{3r,i}(\gamma)}|\Im^{3r,i}(\z)|_{\N}f(|\Im^{3r,i}(\z)|_{\N})>
\frac{t}{3^d.2}\right\rbrack\;.
\end{eqnarray*}
\end{lemma}

\begin{proof}
Let $\gamma$ be any member of $\G$. Recall we say a box is
ugly  if it is not nice, and say it is bad
when it is not good. We write
 $\z\sim \z'$ if $\z$ and $\z'$ are two points adjacent on $\GG_d$, and
$\z\simeq \z'$ if $\z\sim \z'$ or $\z=\z'$. First, we cover $\gamma$ with boxes. This leads to the animal $\A^{r,1}(\gamma)$ and we distinguish two kind of boxes: those all of whose neighbours are nice, and the others. Formally,
$$\Im^{r,1}(\A^{r,1}(\gamma))=\bigcup_{\z\in\A^{r,1}(\gamma)} B_\z^r=U_1\cup U_2\;,$$
where
$$U_1=\bigcup_{\substack{\z\in\A^{r,1}(\gamma)\\\exists
\;\z'\simeq
    \z,\;B_{\z'}^r\mathrm{\;is\;ugly}}} B_\z^r\;,$$
and
$$
U_2=\bigcup_{\substack{\z\in\A^{r,1}(\gamma)\\\forall
\z'\simeq \z,\;B_{\z'}^r\mathrm{\;is\;nice}}} B_\z^r\;.$$
Accordingly,
\begin{eqnarray*}
\nonumber\sum_{v\in \N\cap\gamma}f(d_\N(v))&\leq &S_1+S_2\;.
\end{eqnarray*}
where:
$$S_1=\sum_{v\in\N}f(d_\N(v))\II_{v\in U_1}\;,$$
and
$$S_2=\sum_{v\in\N}f(d_\N(v))\II_{v\in U_2}\;.$$
Now, remark that if there exists $\z'\simeq \z\in \A^{r,1}(\gamma)$ such that
$B_{\z'}^r$ is ugly, then $r\z+B_0^{3r}$ is bad, since it contains some ugly sub-box. Notice that $\{\vec{f}_i+3\ZZ^d,\;i\in\{1,\ldots,3^d\}\}$ is a partition of $\ZZ^d$, so there is a unique pair $(i,\z'')$ in $\{1,\ldots,3^d\}\times\ZZ^d$ such that $\z=3\z''+\vec{f}_i$. Then, $r\z+B_0^{3r}=3r\z''+B_0^{3r}+3r\vec{f}_i/3=\Im^{3r,i}(\z'')$ is a bad box containing $B_\z^r$, and $\z''\in\A^{3r,i}(\gamma)$. This gives:
$$
U_1\subset \bigcup_{i=1}^{3^d}\bigcup_{\substack{\z''\in\A^{3r,i}(\gamma)\\ \Im^{3r,i}(\z'')\mathrm{\;is\;bad}}}  \Im^{3r,i}(\z'')\;.
$$
And thus,
\begin{equation}
\label{eq:S1}
S_1\leq\sum_{i=1}^{3^d}\sum_{\substack{\z\in \A^{3r,i}(\gamma)\\
\Im^{3r,i}(\z)\mathrm{\;is\;bad}}}\sum_{v\in\N\cap
    \Im^{3r,i}(\z)}f(d_\N(v))\;.
\end{equation}
Notice that for $v\in\N$, the degree of $v$ in the Delaunay graph, or equivalently, the number of $(d-1)$-dimensional faces of $C_v$, can be expressed as follows:
\begin{equation}
\label{eq:neighbours}
d_\N(v)=|\{u\in\N\mbox{ s.t. there exists a Delaunay cell containg }u\mbox{ and }v\}|\;.
\end{equation}
When $v$ belongs to some bad box $\Im^{3r,i}(\z)$, Lemma \ref{lemgoodbox} shows that
every Delaunay cell to which $v$ belongs is included in
$\Adh^{3r,i}(\Cl^{3r,i}(\z))$. Therefore, in this case, $d_\N(v)$ is bounded
from above by $|\Adh^{3r,i}(\Cl^{3r,i}(\z))|$. Since $f$ is non-decreasing, we get that for every $(3r,i)$-bad $\z$:
\begin{eqnarray*}
\sum_{v\in\N\cap
    \Im^{3r,i}(\z)}f(d_\N(v))&\leq &\sum_{v\in\N\cap
    \Im^{3r,i}(\z)}f(|\Adh^{3r,i}(\Cl^{3r,i}(\z))|_\N)\;,\\
&=&|\Im^{3r,i}(\z)|_{\N}f(|\Adh^{3r,i}(\Cl^{3r,i}(\z))|_\N)\;.
\end{eqnarray*}
Plugging this into (\ref{eq:S1}) gives:
\begin{eqnarray}
\nonumber S_1&\leq &\sum_{i=1}^{3^d}\sum_{\substack{\z\in \A^{3r,i}(\gamma)\\
\Im^{3r,i}(\z)\mathrm{\;is\;bad}}}|\Im^{3r,i}(\z)|_{\N}f(|\Adh^{3r,i}(\Cl^{3r,i}(\z))|_\N)\;,\\
\nonumber&=&\sum_{i=1}^{3^d}\sum_{\C\in\Cl^{3r,i}(\A^{3r,i}(\gamma))}f(|\Adh^{3r,i}(\C)|_\N)\sum_{\substack{\z\in \A^{3r,i}(\gamma)\\
\Cl^{3r,i}(\z)=\C}}|\Im^{3r,i}(\z)|_{\N}\;,\\
\nonumber&\leq &\sum_{i=1}^{3^d}\sum_{\C\in\Cl^{3r,i}(\A^{3r,i}(\gamma))}f(|\Adh^{3r,i}(\C)|_\N)\sum_{\substack{\z\in \ZZ^d\\
\Cl^{3r,i}(\z)=\C}}|\Im^{3r,i}(\z)|_{\N}\;,\\
\nonumber & = &\sum_{i=1}^{3^d}\sum_{\C\in\Cl^{3r,i}(\A^{3r,i}(\gamma))}f(|\Adh^{3r,i}(\C)|_\N)|\Cl^{3r,i}(\C)|_{\N}\;,\\
\label{eq:sumugly}S_1 &\leq&\sum_{i=1}^{3^d}\sum_{\C\in\Cl^{3r,i}(\A^{3r,i}(\gamma))}|\Adh^{3r,i}(\C)|_{\N}f(|\Adh^{3r,i}(\C)|_\N)\;,
\end{eqnarray}
where we used the fact that $f$ is positive in the last inequality. The last series of inequalities allows us to bound a sum of dependent variables (the first
one) by a sum of variables (the last one) that resembles a greedy animal on $\ZZ^d$.
This will be made clear with Lemma \ref{lemclusterbox}. Now, let us bound $S_2$. If $\z\in\ZZ^d$ is such that for every neighbor $\z'$ of $\z$, $B_{\z'}^r$ is nice, then Lemma \ref{lemtruegoodbox} shows that any Delaunay cell touching $B_{\z}^r$ is included inside $\Adh^{r,1}(\z)$, which is itself included in $\bigcup_{\z'\simeq \z}B_{\z'}^r$. Thus, using equality (\ref{eq:neighbours}):
\begin{eqnarray*}
S_2&\leq& \sum_{\z\in
\A^{r,1}(\gamma)}\sum_{v\in\N\cap
    B_\z^r}f(|\bigcup_{\z'\simeq \z}B_{\z'}^r|_\N)\;,\\
& = & \sum_{\z\in \A^{r,1}(\gamma)}|B_{\z}^r|_\N
f(|\bigcup_{\z'\simeq \z}B_{\z'}^r|_\N)\;,\\
&\leq & \sum_{\z\in \A^{r,1}(\gamma)}|\bigcup_{\z'\simeq \z}B_{\z'}^r|_\N
f(|\bigcup_{\z'\simeq \z}B_{\z'}^r|_\N)\;.
\end{eqnarray*}
Again, $\{\vec{f}_i+3\ZZ^d,\;i\in\{1,\ldots,3^d\}\}$ is a partition of $\ZZ^d$, and for every $\z$ there is a unique couple $(i,\z'')$ in $\{1,\ldots,3^d\}\times\ZZ^d$ such that $\bigcup_{\z'\simeq \z}B_{\z'}^r=\Im^{3r,i}(\z'')$, in which case $\z''\in\A^{3r,i}(\gamma)$. Thus,
\begin{equation}
\label{eq:sumnice}
S_2\leq\sum_{i=1}^{3^d}\sum_{\z\in \A^{3r,i}(\gamma)}|\Im^{3r,i}(\z)|_{\N}f(|\Im^{3r,i}(\z)|_{\N})\;.
\end{equation}
Now, the lemma follows from inequalities (\ref{eq:sumugly}) and (\ref{eq:sumnice}).
\end{proof}

\subsection{The chaining lemma}
\label{subsec:chaininglemma}
The aim of this section is to get good bounds for the
probabilities appearing in the right hand side of the inequality in
Lemma \ref{lemfrontiere}. The following lemma is an adaptation of the technique of \cite{Coxetal93} to control the tail of the greedy animals on $\ZZ^d$.  Notice that Proposition~1 in \cite{Coxetal93} is not enough for us. This stems from the fact that in our setting, according to which weight function $f$ we have, we may hope for better deviation inequalities, simply because the award may have thinner tails than the minimal ones considered in \cite{Coxetal93}. Unfortunately, this necessitates to go over again the whole proof, adjusting the chaining argument. Since this makes the proof quite long and technical, we defer it until section \ref{subsec:animals}.

For each positive integer $m\geq 1$ let:
$$
\Phi_m=\{\A \mbox{ animal in $\GG_d$ s.t.}\,|\A|\leq m,\; \0\in \A\}\;.
$$
\begin{lemma}
\label{lemclusterbox} Let $g$ be a nondecreasing function from $\RR^+$ to
$[1,\infty)$, and define:
$$l(y):=g^{-1}(y)\,,\,\,q(x):=
\hat l(x)\,.
$$
There exists $r_0>0$ such that for all $r>r_0$
there exist positive and finite constants $c_5$, $c_6$ and $c_7$ such that, if:
$$
\liminf_{y\to\infty}\frac{l(y)}{\log y}\geq c_5,
$$
then for every $n\geq m$,
\begin{equation}
\label{eqmajoexpocluster} \PP(\sup_{\A\in\Phi_m}\sum_{\C\in\Cl^{r,i}(\A)}g(|\Adh^{r,i}(\C)|_\N)> c_6n)\leq e^{-c_7l(q(n))}\;,
\end{equation}
and
\begin{equation}
\label{eqmajoexpobox} \PP(\sup_{\A\in\Phi_m}\sum_{\z\in\A}g(|\Im^{r,i}(\z)|_\N)> c_6n)\leq e^{-c_7l(q(n))}\;.
\end{equation}
\end{lemma}


\subsection{Proof of Theorem \ref{majotheo}}
\label{subsec:proofmain}
We shall require the following
result from \cite{Pimentel10} (Theorem 1):
\begin{lemma}
\label{lemmajoanimal-1} There exist finite and positive
constants $z_2$ and $c_{13}$ such that for every $r\geq 1$, and every $i=1,\ldots,3^d$,
$$
\forall z\geq z_2,\;\PP(\max_{\P\in\Pi_n}|\A^{r,i}(\P)|>zn)\leq
e^{-c_{13}zn}\;.
$$
\end{lemma}
Let us prove Theorem \ref{majotheo}. First,
\begin{eqnarray*}
\EE( F_n)=\EE\lbrack \max_{\P\in\Pi_n}\sum_{v\in \P}f(d_\N(v))\rbrack
&=&n\int_0^\infty \PP(\max_{\P\in\Pi_n}\sum_{v\in
\P}f(d_\N(v))>nz)\;dz\;.
\end{eqnarray*}
Let $K$ be a positive real number and $r\geq 1$. We then have 
\begin{eqnarray}
\nonumber\PP(\max_{\P\in\Pi_n}\sum_{v\in \P}f(d_\N(v))>nz)&\leq&\PP(\max_{\P\in\Pi_n}|\A^r_1(\P)|>Knz)\\
\nonumber&+&\PP(\max_{\P\in\Pi_n}\sum_{v\in \P}f(d_\N(v))>nz \mbox{ and }\max_{\P\in\Pi_n}|\A^r_1(\P)|\leq \lfloor Knz\rfloor)\\
\nonumber&\leq&\PP(\max_{\P\in\Pi_n}|\A^r_1(\P)|>Knz)\\
\label{eqmajo-1}&+&\PP(\sup_{\A\in\Phi_{\lfloor
Knz\rfloor}}\sum_{v\in \Im^{r,1}(\A)\cap \N}f(d_\N(v))>nz)\;.
\end{eqnarray}
Thus,
\begin{eqnarray}
\nonumber \EE\lbrack \sum_{v\in\P}f(d_\N(v))\rbrack &\leq &
\frac{1}{K}\EE(\max_{\P\in\Pi_n}|\A^r_1(\P)|)\\
\nonumber &+&n\int_0^\infty\PP(\sup_{\A\in\Phi_{\lfloor
Knz\rfloor}}\sum_{v\in \Im^{r,1}(\A)\cap \N}f(d_\N(v))>nz)\;dz\;,\\
 \nonumber&\leq & \frac{1}{K}\left(\frac{1}{c_{13}}+nz_2\right)\\
  &+&  \label{eqmajokesten} n\int_0^\infty\PP(\sup_{\A\in\Phi_{\lfloor
Knz\rfloor}}\sum_{v\in \Im^{r,1}(\A)\cap \N}f(d_\N(v))>nz)\;dz\;,
\end{eqnarray}
thanks to Lemma \ref{lemmajoanimal-1}.

Now, let $t$ be any positive number and $m$ a positive integer. Using Lemma \ref{lemfrontiere},
\begin{eqnarray*}
 &&\PP(\sup_{\A\in\Phi_m}\sum_{v\in \Im^{r,1}(\A)\cap \N}f(d_\N(v))>t)\\
 &\leq &
\sum_{i=1}^{3^d}\PP\left\lbrack\sup_{\A\in\Phi_m}\sum_{\C\in\Cl^{3r,i}(\A^{3r,i}(\Im^{r,1}(\A)))}|\Adh^{3r,i}(\C)|_{\N}f(|\Adh^{3r,i}(\C)|)>\frac{t}{3^d.2} \right\rbrack\\
&&+\sum_{i=1}^{3^d}\PP\left\lbrack\sup_{\A\in\Phi_m}\sum_{\z\in \A^{3r,i}(\Im^{r,1}(\A))}|\Im^{3r,i}(\z)|_{\N}f(|\Im^{3r,i}(\z)|_{\N})>
\frac{t}{3^d.2}\right\rbrack\;,
\end{eqnarray*}
Notice that for each $\A$ in $\Phi_m$, $\A^{3r,i}(\Im^{r,1}(\A))$ will be an animal containing $\0$ and with less vertices than $\A$. Thus, there is some $\A'$ in $\Phi_m$ such that $\A^{3r,i}(\Im^{r,1}(\A))\subset \A'$, which yields to
\begin{eqnarray}
\nonumber &&\PP(\sup_{\A\in\Phi_m}\sum_{v\in \Im^{r,1}(\A)\cap \N}f(d_\N(v))>t)\\
\nonumber &\leq &
\sum_{i=1}^{3^d}\PP\left\lbrack\sup_{\A'\in\Phi_m}\sum_{\C\in\Cl^{3r,i}(\A')}|\Adh^{3r,i}(\C)|_{\N}f(|\Adh^{3r,i}(\C)|)>\frac{t}{3^d.2} \right\rbrack\\
\label{eqclusterbox}&&+\sum_{i=1}^{3^d}\PP\left\lbrack\sup_{\A'\in\Phi_m}\sum_{\z\in \A'}|\Im^{3r,i}(\z)|_{\N}f(|\Im^{3r,i}(\z)|_{\N})>
\frac{t}{3^d.2}\right\rbrack\,.
\end{eqnarray}
Now, if
$$
\liminf_{y\to\infty}\frac{\hat f(y)}{\log y}\geq c_5\;,
$$
then inequality (\ref{eqclusterbox}) and Lemma \ref{lemclusterbox}
imply:
$$
\PP(\sup_{\A\in\Phi_m}\sum_{v\in \Im^{r,1}(\A)\cap \N}f(d_\N(v))>c_6m)\leq 3^d.2e^{-c_7f^*(m)}\;.
$$
Therefore, fix $K=\frac{1}{c_6}$, 
$c_1=\sup\{\frac{2}{c_7},c_5\}$ and suppose that
$$
\liminf_{y\to\infty}\frac{\hat f(y)}{\log y}\geq c_1\;.
$$
This implies that $e^{-c_7f^*(x)}$ is integrable on $\RR^+$. Thus,
\begin{eqnarray*}
\PP(\sup_{\A\in\Phi_{\lfloor
Knz\rfloor}}\sum_{v\in \Im^{r,1}(\A)\cap \N}f(d_\N(v))>nz)&\leq&3^d.2e^{-c_7f^*(Knz)}\;,
\end{eqnarray*}
and then,
$$
\int_0^\infty\PP(\sup_{\A\in\Phi_{\lfloor
    Knz\rfloor}}\sum_{v\in\Im^{r,1}(\A)\cap \N}f(d_\N(v))>nz)\;dz=O\left(\frac{1}{n}\right)\;.
  $$
Plugging these bounds into (\ref{eqmajo-1}), (\ref{eqmajokesten})
leads to the desired result. \hfill $\square$

\subsection{Sketch of the proof of Proposition \ref{prop:Gamma}}
\label{subsec:proofgamma}

It turns out that the proof of Proposition \ref{prop:Gamma} is very similar
to that of Theorem \ref{majotheo}, so we shall only stress the
differences. First, we need to see that the notion of a nice box is still a good one to
control $\Gamma(\N,\phi)$. To do that, we need to express
$\Gamma(\N,\phi)$ differently. Let $\E(\phi)$ be the set of edges of
$\phi$. Then,
$$\Gamma(\N,\phi)=\bigcup_{e\in\E(\phi)}\Gamma(\N,e)\;,$$
where:
$$\Gamma(\N,e)=\{x\in\RR^d\mbox{ s.t. }e\not\subset
\mathcal{D}(\N\cup \{x\}\}\;.$$
Recall that $B(\Delta)$ denotes the closed circumball of a Delaunay cell $\Delta$. It may be seen that (cf. Lemma \ref{lemm:gammaNe}):
$$\Gamma(\N,e)=\bigcap_{\Delta\ni e}\inter{B}(\Delta)\;,$$
where the intersection is taken  over all Delaunay cells $\Delta$
which contain $e$. We get easily the following result from Lemma \ref{lemtruegoodbox}.
\begin{lemma}
\label{lemgoodboxGamma}
Fix $r>0$ and $i\in\{0,\ldots,3^d\}$, assume that $\C$ is an $(r,i)$-ugly cluster in $\GG_d$. Let $\phi$ be a self-avoiding
path in $\D(\N)$. Then, for any vertex $v$ in $\phi$, and any edge
$e\in\E(\phi)$ such that $v\in e$, 
$$v\in \Im^{r,i}(\C)\Rightarrow \Gamma(\N,e)\subset
\Adh^{r,i}(\C)\;.$$
The same holds for bad clusters.
\end{lemma}
\begin{proof}
We may drop $r$ and $i$ for simplicity. Let $B'$ be any circumball of a Delaunay cell containing $v$. Then, $\C$ and $B'$ satisfy the hypotheses of Lemma \ref{lemtruegoodbox}, and thus $B'\subset  \Adh(\C)$. Thus, for any edge
$e\in\E(\phi)$ such that $v\in e$, we have $\Gamma(\N,e)\subset  \Adh(\C)$.
\end{proof}
Notice that a vertex in $\GG_d$ has $3^d-1$ neighbours. Thus, the number of points in the union of a cluster $C$ and its exterior boundary is at most $3^d|C|$. 
Lemma \ref{lemgoodboxGamma} shows that when $v\in \Im^{3r,i}(\C)$,
then, $\mu(\Gamma(\N,e))\leq \mu(\Adh^{3r,i}(\C))\leq c_\mu(3r)^d3^d|\C|$. So
we obtain the following analogue of Lemma \ref{lemfrontiere}.
\begin{lemma}
\label{lemfrontiereGamma} 
\begin{eqnarray*}
\PP\lbrack\max_{\phi\in \SA_n}\mu(\Gamma(\N,\phi))>
t\rbrack&\leq&
\sum_{i=1}^{3^d}\PP\left\lbrack\max_{\phi\in \SA_n}\sum_{\C\in\Cl^{3r,i}(\A^{3r,i}(\phi))}|\C|>\frac{t}{3^d.2.c_\mu(9r)^d} \right\rbrack\\
&+&\sum_{i=1}^{3^d}\PP\left\lbrack\max_{\phi\in \SA_n}|\A^{3r,i}(\phi)|>
\frac{t}{3^d.2.c_\mu(9r)^d}\right\rbrack\;.
\end{eqnarray*}
\end{lemma}
The rest of the proof is similar to that of Theorem \ref{majotheo}. Notice first that $\SA_n\subset\Pi_n$. Lemma \ref{lemmajoanimal-1} implies that for all $r\geq 1$ and $z\geq z_2$,
\begin{equation}
\label{eqA}
\sum_{i=1}^{3^d}\PP\left\lbrack\max_{\phi\in \SA_n}|\A^{3r,i}(\phi)|>
zn\right\rbrack\leq 3^de^{-c_{13}zn}\;.
\end{equation}
Now we fix $r=r_0$ large enough so that we are in the sub-critical phase for the percolation of clusters of bad boxes (see the beginning of section \ref{subsec:animals} for more details). Then, there is some $\lambda>0$ such that:
$$ \EE[e^{\lambda|\Cl(0)|}]=:M<\infty\;.$$
Lemma \ref{lemmajoanimal-1} shows that for any $z\geq z_2$ and any $i$:
\begin{equation}
\label{eqCAC}
\PP\left\lbrack\max_{\phi\in \SA_n}\sum_{\C\in\Cl^{3r,i}(\A^{3r,i}(\phi))}|\C|>zn \right\rbrack \leq  e^{-c_{13}zn}+\PP\left\lbrack\max_{\Lambda\in \Phi_{\lfloor nz\rfloor}}\sum_{\C\in\Cl(\Lambda)}|\C|>zn \right\rbrack\;
\end{equation}
Now, using Lemma \ref{lemcorrelneg}, one gets
\begin{eqnarray*}
\PP\left\lbrack\max_{\Lambda\in \Phi_{\lfloor nz\rfloor}}\sum_{\C\in\Cl(\Lambda)}|\C|>zn \right\rbrack &\leq &\sum_{\Lambda\in \Phi_{\lfloor nz\rfloor}}\PP\left\lbrack\sum_{\C\in\Cl(\Lambda)}|\C|>zn \right\rbrack\;,\\
&\leq & \sum_{\Lambda\in \Phi_{\lfloor nz\rfloor}}e^{-\lambda nz}\EE\left\lbrack e^{\lambda\sum_{\C\in\Cl(\Lambda)}|\C|}\right\rbrack\;,\\
&\leq & \sum_{\Lambda\in \Phi_{\lfloor nz\rfloor}}e^{-\lambda nz}M^{|\Lambda|}\;.
\end{eqnarray*}
It is well known that there is a finite constant $K>1$, depending on the dimension such that $| \Phi_k|\leq K^k$ for every $k$. Thus:
$$\PP\left\lbrack\max_{\Lambda\in \Phi_{\lfloor nz\rfloor}}\sum_{\C\in\Cl(\Lambda)}|\C|>zn \right\rbrack \leq (KM)^{nz}e^{-\lambda nz}\;.$$
Thus, for any $z$ larger that $2\log(KM)/\lambda$,
\begin{equation}
\label{eqPhiC}
\PP\left\lbrack\max_{\Lambda\in \Phi_{\lfloor nz\rfloor}}\sum_{\C\in\Cl(\Lambda)}|\C|>zn \right\rbrack \leq e^{-\lambda nz/2}\;.
\end{equation}
Gathering equations (\ref{eqA}), (\ref{eqCAC}) and (\ref{eqPhiC}), we conlude that for $z$ large enough, and for $c=\min\{c_{13}/(3^d.2.c_\mu(9r)^d),\lambda/2\}$, we have for any $n$:
$$\PP\lbrack\max_{\phi\in \SA_n}\mu(\Gamma(\N,\phi))> t\rbrack\leq 3^{d+1}e^{-cnz}\;.$$

\section{Applications}
\label{sec:applications}
\subsection{The connectivity constant on the Delaunay triangulation}
\label{connectivity}
Problems related to self-avoiding paths are connected with various
branches of applied mathematics such as long chain polymers,
percolation and ferromagnetism. One fundamental problem is the
asymptotic behavior of the connective function $\kappa_n$, defined by
the logarithm of the number $N_n$ of self-avoiding paths (on a
fixed graph $\G$) starting at vertex $v$ and with $n$ steps. For
planar and periodic graphs subadditivity arguments yields that
$n^{-1}\kappa_n$ converges, when $n\to\infty$, to some value
$\kappa\in(0,\infty)$ (the connectivity constant) independent of the
initial vertex $v$. In disordered planar graphs subadditivity is lost
but, if the underline graph possess some statistical symmetries
(ergodicity), we may believe that the rescaled connective function
still converges to some constant. When $x$ is in $\RR^d$, let us denote
by $v(x)$ the (a.s unique) point of $\N$ such that $x\in C_v$. Let $N_n$ denote the number of
 self-avoiding paths of length $n$ starting at $v(0)$ on the Delaunay
triangulation of a Poisson random set $\N$. We recall that the
intensity of the Poisson random set is bounded from above and below by a constant times the Lebesgue measure on $\RR^d$. From Theorem \ref{majotheo} we obtain a linear upper bound for the connective function of the Delaunay triangulation:

\begin{proposition}
\label{prop:connectivity}
Let $\kappa_n=\log N_n$. There exist positive constants $z_3$ and $c_2$ such that, for every
$u\geq nz_3$,
$$\PP(\kappa_n\geq u)\leq e^{-c_2 u^{1/3}}\;.$$
In particular,
$$ \EE(\kappa_n)=O(n)\;.$$
\end{proposition}
\begin{proof}
We shall use the following intuitive inequality which is a consequence
of Lemma \ref{lemm:connect} below.
$$N_n\leq \sup_{\P\in\Pi_{n-1}}\prod_{v\in \N\cap \P}d_\N(v)\;.$$
We now deduce the following:
\begin{eqnarray*}
\PP(N_n\geq t)&\leq &\PP(\sup_{\P\in\Pi_{n-1}}\prod_{v\in \N\cap
  \P}d_\N(v)\geq t)\;,\\
&=&\PP(\sup_{\P\in\Pi_{n-1}}\sum_{v\in \N\cap
  \P}\log (d_\N(v))\geq \log t)\;,\\
&\leq&\PP(\sup_{\P\in\Pi_{n-1}}\sum_{v\in \N\cap
  \P} d_\N(v)\geq \log t)\;.
\end{eqnarray*}
Now, notice that $f(x)=x$ is $c-$nice for every $c>0$. Thus, according
to Theorem \ref{majotheo}, there are positive constants $z_3$ and $c_2$ such that, for all $t\geq e^{rz_3}$,
$$\PP(N_n\geq t)\leq e^{-c_2 (\log t)^{1/3}}\;.$$
Or, equivalently, for every $u\geq rz_3$,
$$\PP(\kappa_n\geq u)\leq e^{-c_2 u^{1/3}}\;.$$
This implies notably that:
$$\EE(\kappa_n)=O(n)\;.$$
\end{proof}
\begin{lemma}
\label{lemm:connect}
Let $G=(V,E)$ be a graph with set of vertices $V$ and set of edges
$E$. Define, for any $n\in\NN$, any $x\in V$ and $I\subset V$:
$$\Delta_n(x,I)=\{\mbox{ s.a paths }\gamma\mbox{ with }n\mbox{ edges
  and s.t. }\gamma_0=x,\;\gamma\cap I=\emptyset\}\;.$$
Define:
$$N_n(x,I)=|\Delta_n(x,I)|\;.$$
Then, for any $n\geq 1$,
$$\forall x\in V,\;\forall I\subset V,\;N_n(x,I)\leq
\sup_{\gamma\in\Delta_{n-1}(x,I)}\prod_{v\in\gamma}d(v)\;,$$
where $d(v)$ stands for the degree of the vertex $v$.
\end{lemma}
\begin{proof}
We shall prove the result by induction. When $n=1$ it is obviously
true. Indeed, if $x\in I$, then $N_1(x,I)=0$ and if $x\not\in I$,
$N_1(x,I)\leq d(x)$. Suppose now that the result is true until
$n\geq 1$. Let $x\in V$ and $I\subset V$. If $x\in I$, then
$N_{n+1}(x,I)=0$. Suppose thus that $x\not\in I$. Then, denoting $u\sim x$ the fact that $u$ and $x$ are
neighbours, and using the induction hypothesis,
\begin{eqnarray*}
N_{n+1}(x,I)&=&\sum_{u\sim x}N_n(u,I\cup\{x\})\;,\\
&\leq&\sum_{u\sim x}\sup_{\gamma\in\Delta_{n-1}(u,I\cup\{x\})}\prod_{v\in\gamma}d(v)\;,\\
&\leq&\sum_{u\sim x}\sup_{u\sim
  x}\sup_{\gamma\in\Delta_{n-1}(u,I\cup\{x\})}\prod_{v\in\gamma}d(v)\;,\\
&=&\sup_{u\sim
  x}\sup_{\gamma\in\Delta_{n-1}(u,I\cup\{x\})}d(x)\prod_{v\in\gamma}d(v)\;,\\
&=&\sup_{\gamma'\in\Delta_{n}(x,I)}\prod_{v\in\gamma'}d(v)\;,
\end{eqnarray*}
and the induction is proved.
\end{proof}
\subsection{First passage percolation on the Poisson-Delaunay graph}
\label{subsec:FPP}
On any graph $\G$, one can define a First Passage Percolation model
(FPP in the sequel)
as follows. Assign to each edge $e$ of $\G$ a random non-negative time
$t(e)$ necessary for a particle to cross edge $e$. The First Passage
time from a vertex $u$ to a vertex $v$ on $\G$ is the minimal time 
needed for a particle to go from $u$ to $v$ following a path on
$\G$. This is of course a random time, and the understanding of the
typical fluctuations of this times when $u$ and $v$ are far apart is
of fundamental importance, notably because of its physical interpretation as a
growth model. The model was first introduced by \cite{HammersleyWelsh65} on $\ZZ^d$ with i.i.d passage times. Some variations on this
model have been proposed and studied, see \cite{Howard04} for a review. One
of these variations is FPP on the Delaunay triangulation, introduced
in \cite{VahidiAslWierman87}, where the
graph $\G$ is the Delaunay graph of a Lebesgue-homogeneous Poisson
random set in $\RR^d$. Often, this Delaunay triangulation is
heuristically considered to behave like a perturbation of the triangular lattice. One important question is whether the fact that
$\G$ itself is random affects the fluctuations of the first passage
times. Since it is not already known exactly what is the right order
of these fluctuations on the triangular lattice, one does not have an
exact benchmark, but the best bound known up today is that the
variance of the passage time between the origin and a vertex at
distance $n$ is of order $O(n/\log n)$. We are not able to show a
similar bound for FPP on the Delaunay triangulation, but we can show
the analogue of Kesten's bound of \cite{Kesten93}: the above variance is at
most of order $O(n)$. This improves upon the results in \cite{Pimentel05} and
answers positively to open problem~12, p.~169 in \cite{Howard04}. To prove this result, we shall also need to bound from above the length of the optimal path, as was done in Proposition~5.8 in \cite{Kesten84} for the classical model. This requires a separate argument, given in section \ref{subsubsec:linearlength}. The bound on the variance will be given in section \ref{subsubsec:varianceFPP}.

\subsubsection{Linear passage weights of self-avoiding paths in percolation}
\label{subsubsec:linearlength}

When $\D$ is a fixed Delaunay triangulation on the plane, the bond percolation model $X=\{X_e\,:\,e\in\D\}$ on the Delaunay triangulation, with parameter $p$ and probability law $\PP_p(.|\D)$, is constructed by choosing each edge $e\in\D$ to be
open, or equivalently $X_e=1$, independently with probability $p$ and
closed, equivalently $X_e=0$, otherwise. An open path is a path
composed by open edges. We denote by $\PP_p$ the probability measure
obtained when $\D$ is the Delaunay triangulation based on $\N$, which is distributed like a Poisson random set on $\RR^d$ with
intensity measure $\mu$. We denote by $\calC_\0$ the maximal connected
sugraph of $\D$ composed by edges $e$ that belong to some open path
starting from $v_\0$, the point of $\N$ whose Voronoi tile contains the origin, and we call it the open cluster. Define the critical probability:
\begin{equation}\label{crit}
\bar{p}_c(d)=\sup\Big\{p\in[0,1]\,:\,\forall \alpha >0\sum_{n\geq 1}n^{\alpha}\PP_p(|\calC_\0|=n)<\infty\Big\}\,.
\end{equation}
It is extremely plausible that this critical probability coincides
with the classical one, i.e $\sup\{p\mbox{
  s.t. }\PP_p(|\calC_\0|=\infty)= 0\}$, but we are not going to prove
this here. By using Proposition \ref{prop:connectivity} to count the
number of self-avoiding paths of size $n$, and then evaluating the
probability that it is an open path,  there exists a
positive and finite constant $B=B(d)$ such that if $p<1/B$ then the
probability of the event $\{|\calC_\0|=n\}$ decays like $e^{-\alpha
  n^{1/(3d)}}$ for some constant $\alpha$. Consequently: 
\begin{lemma}\label{lcrit}
$0<\frac{1}{B}\leq \bar{p}_c(d)$.
\end{lemma} 
The main result of this section is the following estimate on the density of open edges.
\begin{theorem}
\label{theo:lden}
If $(1-p)<\bar{p}_c(d)$ then there exists constants (only depending on $p$) $a_1,a_2>0$ such that  
$$\PP_p\Big(\exists \mbox{ s.a. }\gamma\mbox{ s.t. }|\gamma|\geq m\mbox{ and
}\sum_{e\in\gamma}X_e\leq a_1m\Big)\leq e^{-a_2m}\,.$$
\end{theorem}

\begin{remark} 
We could prove Theorem \ref{theo:lden} directly under $(1-p)<1/B$ (again using Proposition \ref{prop:connectivity}) but we want to push optimality as far as we can.
\end{remark}

To prove Theorem \ref{theo:lden} we shall first obtain a control on the
number of boxes intersected by a self-avoiding path in $\D$. Recall
that, in section \ref{comp}, for fixed $L>0$ and for each self-avoiding path $\gamma$ in $\D$ we have defined an animal $\A^r(\gamma):=\A^{r,1}(\gamma)$ on $\GG^d$ by taking the vertices $\z$ such that $\gamma$ intersects $B_\z^{r}$. By Corollary 2 in \cite{Pimentel10} we have: 
\begin{lemma}\label{l5}
For each $r\geq 1$ there exist finite and positive constants $b_3,b_4,b_5,b_6$ such that for any $x\geq 0$ (check reference)
\begin{equation}\label{e1l5}
\PP\left(\min_{\v_0\in\gamma\,,|\gamma|\geq r\,}|\A^{r}(\gamma)|<b_3 x\right)\leq e^{-b_4 x}\, ,
\end{equation}
and
\begin{equation}\label{e0l5}
\PP\left(\max_{\v_0\in\gamma\,,|\gamma|\leq r\,}|\A^{r}(\gamma)| >b_5 x\right)\leq e^{-b_6 x}\, .
\end{equation}
\end{lemma}

{\it Proof of Theorem \ref{theo:lden}:~} Let $L>0$, $\z\in\GG^d$. In this
section, we call $B_\z^{L/2}$ a \emph{good box} if:
\begin{itemize}
\item[(i)] For all $\z'\in\GG^d$ with $|\z-\z'|_\infty = 2$,
  $B_{\z'}^{L/2}$ is $\N$-nice (recall Definition \ref{defi:goodbox}),
\item[(ii)] For all $\gamma$ in $\D$ connecting $B_\z^{L/2}$ to $\partial B_\z^{3L/2}$ we have $\sum_{e\in\gamma}X_e \geq 1$\,.
\end{itemize} 
Let 
$$Y_\z(L)=\II_{B_\z^{L/2}\mbox{ is a good box }}.$$
Our first claim is: if $1-p< \bar{p}_c(d)$ then
\begin{equation}\label{l1}
\lim_{L\to\infty}\PP_p\big( Y_\z(L)=1\big)=1\,.
\end{equation}
Since the intensity of the underlying Poisson random set is comparable
with the Lebesgue measure (\ref{eq:comparaison}), condition (i) has
probability going to one as $L$ goes to infinity. Now, denote by $E_L$
the event that (ii) is false and (i) is true. Suppose that $E_L$
occurs, and let $\gamma$ be a
path contradicting (ii). We may choose the first edge of $\gamma$ as
some $[v_1,v_2]$ intersecting $B_\z^{3L/2}$. Thanks to Lemma
\ref{lemgoodbox}, and since (i) is true, $[v_1,v_2]$ lies in
$B_\z^{5L/2}$. Also, there is a point $v'$ such that either $v_1$ or
$v_2$ lies at distance at least $L$ from $v'$ and $\gamma$ connects
$v_1$ and $v_2$ to $v'$ with closed edges. Divide $B_\z^{5L/2}$
regularly into subcubes of side length 1. This partition $\P_L$ has
cardinality of order $L^d$. For each box $B$, let $A(B)$ be
the event that there exist a vertex $v\in \N\cap B$, $v'\in \N$ such
that $|v-v'|\geq L$, and a path $\gamma$ in $\D$ from $v$ to $v'$ such
that $\sum_{e\in\gamma}X_e=0$. Denote by $A_L$ the event
$A(B_\0^{1/2})$. The remarks we just made imply that:
$$\PP_p\big( E_L\big)\leq \sum_{B\in\P_L}\PP_p\big(A(B)\big)\leq cL^{d} \PP_p\big(A_L\big)\,.$$
Now, let us bound $\PP_p\big(A_L\big)$. If $p=1$,
$\PP_p\big(A_L\big)=0$, so we suppose that $p<1$. Let $\B_d$ denote the ball of
center $\0$ and radius $3\sqrt{d}$, and let $\E_d$ be the set of edges
of $\D$ lying completely inside $\B_d$. A simple geometric lemma (see
Lemma \ref{lem:cheminborne} in section \ref{subsec:lemchemborne}) shows that if $v\in N\cap B_0^{1/2}$,
there is a path from $v_\0$ to $v$ on $\D$ whose edges are in
$\E_d$. We can further restrict this path to belong to a spanning tree $\mathcal{T}_d$ of the connected component of $v_\0$ in $\E_d$. Thus, if $A_L$ holds and if all the edges of $\mathcal{T}_d$ are closed,
then there is a closed path from $v_\0$ to some vertex $v'$ such that
$|v'|\geq L-\sqrt{2}$. Let $L$ be larger than $2\sqrt{2}$. Using
(\ref{e0l5}) (with $l=1$ and $r$ a constant times $L$), there exist constants $a,b>$ such that,
with probability greater than $(1-e^{-aL})$, a path $\gamma$ in  the
Poisson-Delaunay triangulation, that connects $v_\0$ to a point $v'$
such that $|v'|\geq L-\sqrt{2}$, has at least $bL$ edges. Now, using
the FKG inequality, we have:
$$\PP_p(A_L|\D)\leq\PP_p(A_L|\{\forall e\in\D,\;X_e=0\},\;\D)\leq
\frac{1}{(1-p)^{|\mathcal{T}_d|}}\left(\PP_p(|\C^{cl}_\0|\geq
    bL|\D)+e^{-aL}\right)\;,$$
where $\C^{cl}_\0$ is the \emph{closed} cluster containing $\0$. Thus:
$$\PP_p(A_L)\leq \EE\left[(1-p)^{-|\mathcal{T}_d|}\left(\PP_p(|\C^{cl}_\0|\geq
    bL|\D)+e^{-aL}\right)\right]\;.$$
Notice that for any $x>0$,
\begin{eqnarray*}
\EE\left[(1-p)^{-|\mathcal{T}_d|}\PP_p(|\C^{cl}_\0|\geq bL|\D)\right]&\leq& \EE\left[
(1-p)^{-|\mathcal{T}_d|}\II_{|\mathcal{T}_d|\geq
  x}\right]\;,\\
&&+\EE\left[
(1-p)^{-|\mathcal{T}_d|}\II_{|\mathcal{T}_d|<
  x}\PP_p(|\C^{cl}_\0|\geq
    bL|\D)\right]\;,\\
&\leq &\EE\left[(1-p)^{-|\mathcal{T}_d|}\II_{|\mathcal{T}_d|\geq
 x}\right]+(1-p)^{-x}\PP_p(|\C^{cl}_\0|\geq bL)\;.
\end{eqnarray*}
 Now, $|\mathcal{T}_d|$ is less than $|\B_d|_{\N}$ and notice that $|\B_d|_{\N}$ has Poisson
distribution with parameter $\mu(\B_d)$. Thus, one can find some
positive constants $C_1$, $C_2$ and $C_3$, depending only on $\mu$ and
$d$, such that, for every $L>0$:
\begin{equation}
\label{eq:controleNB_d}
\PP(|\mathcal{T}_d|\geq C_1\log(L))\leq \PP(|\B_d|_{\N}\geq C_1\log(L))\leq C_2L^{-2(d+1)}\;.
\end{equation}
Now, using Cauchy-Schwartz inequality, and inequality (\ref{eq:controleNB_d}),
\begin{eqnarray*}
\PP_p(A_L)&\leq& e^{-aL}\EE\left[(1-p)^{-|\mathcal{T}_d|}\right]\\
&&+\sqrt{\EE\left[(1-p)^{-2|\mathcal{T}_d|}\right]}.\sqrt{\PP(|\mathcal{T}_d|\geq
  C_1\log(L))}+L^{C_1\log\frac{1}{1-p}}\PP_p(|\C^{cl}_\0|\geq bL)\;,\\
&\leq& e^{-aL}C_3+C_3\sqrt{C_2}L^{-(d+1)}+L^{C_1\log\frac{1}{1-p}}\PP_p(|\C^{cl}_\0|\geq bL)\;,
\end{eqnarray*}
where $C_3$ is some finite bound (depending on $p<1$) on
$\sqrt{\EE\left[(1-p)^{-2|\mathcal{T}_d|}\right]}$ (all exponential moments of a Poisson random variable are finite). Now, if $0<(1-p)<\bar{p}_c(d)$,
\begin{eqnarray}
\nonumber L^{d}\PP\big(A_L\big)&\leq& L^de^{-aL}C_3+C_3\sqrt{C_2}L^{-1}+L^{d+C_1\log\frac{1}{1-p}}\PP_p(|\C^{cl}_\0|\geq bL)\;,\\
\nonumber &\leq& \sum_{n\geq bL}\left(\frac{n}{b}\right)^{d+C_1\log\frac{1}{1-p}}\PP_p(|\C^{cl}_\0|=n)+L^de^{-aL}C_3+C_3\sqrt{C_2}L^{-1}\;,
\end{eqnarray}
which goes to zero as $L$ goes to infinity. This concludes the proof of (\ref{l1}).

Our second step is: the collection $\{Y_\z(L)\,:\,\z\in\GG^d\}$ is a $5$-dependent Bernoulli field. To see this, first notice that $Y_z(L)= Y_\z(L,\N,X)$. Thus, the claim will follow if we prove that $Y_\z(L,\N,X)=Y_\z(L,\N\cap B_\z^{5L/2},X)$. To prove this, assume that $Y_\z(L,\N,X)=0$. Then,
either (i) does not hold, or (i) holds and (ii) does not. To see if (i) does not hold, it is clear
that we only have to check $\N\cap B^{5L/2}_\z$. If (i) holds and (ii) not, we then use Lemma 2.2 to show that, under (i), inside $B^{3L/2}_\z$ the Delaunay triangulation based on $\N$ is the same as the one based on $\N\cap B^{5L/2}_\z$. Thus, if (ii) does not hold for $\N$, it will certainly not hold for $\N\cap B^{5L/2}_\z$. The same argument works to show that if $Y_\z(L,\N,X)=1$ then $Y_\z(L,\N \cap B^{5L/2}_\z,X) = 1$.

With (\ref{l1}) and $5$-dependence in hands, we can chose $L_0\geq 1$ large enough such that 
\begin{equation}\label{e1}
\PP_p\left(\min_{\A\in\Phi^r}\sum_{\z\in\A}Y_\z<c_1r\right)\leq e^{-c_2 r}\,,
\end{equation}
for $c_1,c_2$ only depending on $L_0$, where $\Phi^r$ denotes the set of lattice animals of size $\geq r$ (see Lemma 9 of \cite{Pimentel10}). 

Let $\gamma$ be a path in $\D$ with $v_0\in\gamma$ and $|\gamma|\geq
m$. Let 
$$\mathcal{S}_L(\gamma):=\{\z\in\A^L(\gamma)\,:\,Y_\z(L)=1\}\,.$$
Notice that there exists at least one subset $\mathcal{S}$ of $\mathcal{S}_L(\gamma)$ such that $|\z-\z'|_\infty\geq 5$ for all $\z,\z'\in\mathcal{S}$ and $k=|\mathcal{S}|\geq |\mathcal{S}_L(\gamma)|/5^d$. Now, write $\mathcal{S}=\{\z_1,\dots,\z_k\}$. By Lemma \ref{lemgoodbox}, one can find disjoint pieces of $\gamma$, say $\gamma_1,\dots,\gamma_k$ such that, for $i=1,\dots,k$, $\sum_{e\in\gamma_i}X_e\geq 1$. Hence,
$$\sum_{e\in\gamma} X_e\geq \sum_{i=1}^k\left(\sum_{e\in\gamma_i}X_e\right)\geq |\mathcal{S}|\geq  \frac{|\mathcal{S}_L(\gamma)|}{5^d}= \frac{\sum_{\z\in\A^L(\gamma)}Y_z(L)}{5^d}\,,$$
which shows that 
\begin{eqnarray}
\nonumber\PP\left(\exists \mbox{ s.a. }\gamma\mbox{ s.t. }|\gamma|\geq m\mbox{ and
}\sum_{e\in\gamma}X_e\leq a_1m\right)&\leq& \PP\left( \min_{\v_0\in\gamma\,,|\gamma|\geq m\,}|\A^L(\gamma)|<b_3m\right)\\
\nonumber&+&\PP\left(\min_{\A\in\Phi^{b_3m}} \sum_{\z\in\A}Y_\z\leq5^d a_1m \right)\,.
\end{eqnarray}
Combining this together with (\ref{e1l5})
and (\ref{e1}), we finish the proof of Theorem \ref{theo:lden}. 

{\hfill $\square$}

\subsubsection{The variance of the first passage time is at most linear} 
\label{subsubsec:varianceFPP}
 Let $\nu$ be
a probability measure on $\RR^+$. Recall that $\Omega$ is the set of locally finite subsets of $\RR^d$. Let $\pi$ denote the Poisson measure on
$\RR^d$ with intensity $\mu$, where $\mu$ is comparable to the
Lebesgue measure on $\RR^d$, in the sense of (\ref{eq:comparaison}). Let $\ER$ denote the set of pairs $\{x,y\}$ of points of
$\RR^d$. We endow the space $\Omega\times\RR_+^{\ER}$ with the product measure $\pi\otimes
\nu^{\otimes \ER}$ and each element of $\Omega\times\RR_+^{\ER}$ is
denoted $(\N,\tau)$. This means that $\N$ is a Poisson random set with intensity $\mu$, and to each edge $e\in \D(\N)$ is
independently assigned a non-negative random variable $\tau(e)$ from the common
probability measure $\nu$.

The passage time $T(\gamma)$ of a path $\gamma$ in the Delaunay triangulation is the sum of the passage times of the edges in $\gamma$:
$$
T(\gamma)=\sum_{e\in\gamma}\tau(e)\;.$$ 
The first-passage time between two vertices $v$ and $v'$ is defined by
$$T(v,v'):=T(v,v',\N,\tau))=\inf\{T(\gamma)\,;\,\gamma\in\Gamma(v,v',\N)\}\;.$$
where $\Gamma(v,v',\N)$ is the set of all finite paths connecting $v$
to $v'$. Given $x,y\in\RR$ we define $T(x,y):=T(v(x) ,v(y))$. 

Remark that $(\N,\tau)\mapsto T(x,y)$ is measurable with respect
to the completion of $\mathcal{F}\otimes \mathcal{B}(\RR_+)^{\otimes\ER}$, where
$\mathcal{B}(\RR_+)$ denotes the Borel $\sigma$-field over $\RR_+$ and
$\mathcal{F}$ is the smallest algebra on $\Omega$ which lets the
coordinate applications $\N\mapsto \II_{x\in\N}$
measurable. To see this, fix $x,y\in\RR^d$ and for each $r>|x,y|$ define $T_r(x,y)$ to be the first-passage time restricted to the Delaunay graph $\D(\N\cap D_r(x))$, where $D_r(x)$ is the ball centred at $x$ and of radius $r$. Then clearly $T_r$ is measurable, since it is the infimum over a countable collection of measurable functions. On the other hand, $T_r(x,y)$ is non-increasing with $r$, and so $T=\lim_{r\to\infty}T_r(x,y)$ is also measurable.
   
We may now state the announced upper-bound on the typical fluctuations
of the first passage time.  We denote by $M_k(\nu)$ the $k-th$ moment of $\nu$:
$$M_k(\nu)=\int |x|^k\;d\nu(x)\;.$$
\begin{theorem}
\label{thmmajovariance}
For any natural number $n\geq 1$, let $\vec{n}$ denote the point
$(n,0,\ldots,0)\in\RR^d$. Assume that $\nu(\{0\})<\bar{p}_c$ and that 
$$M_2(\nu):=\int x^2\;d\nu(x)<\infty\;.$$ Then, as $n$ tends to infinity,
$$\Var(T(0,\vec{n},\tau,\N))=O(n)\;.$$
\end{theorem}

Before we prove Theorem \ref{thmmajovariance}, we need to state a Poincar\'e
inequality for Poisson random sets (for a proof, see for instance
\cite{HoudrePrivault03}, inequality (2.12) and Lemma 2.3).

\begin{proposition}
\label{prop:PoincarePoisson}
Let $F:\Omega\rightarrow \RR$ be a square integrable random
variable. Then,
$$\Var_{\pi}(F)\leq\EE_{\pi}(V_-(F))\;,$$
where:
$$V_-(F)(\N)=\sum_{v\in \N}(F(\N)-F(\N\setminus\{v\}))_-^2+\int
(F(\N)-F(\N\cup \{x\}))_-^2\;d\lambda(x)\;.$$
\end{proposition}  

A geodesic between $x$ and $y$, in the FPP model, is a path $\rho(x,y)$ connecting $v_x$ to $v_y$ and such that 
$$
T(x,y)=T\big(\rho(x,y)\big)=\sum_{e\in\rho(x,y)}\tau(e)\,.
$$

When a geodesic exists, we shall denote by $\rho_n(\tau,\N)$ a geodesic between $0$ and $\vec{n}$ with
minimal number of vertices. If there are more than one such geodesics,
we select one according to some deterministic rule. We shall write $T_n(\tau,\N)$ for
$T(0,\vec{n},\tau,\N)$.

To control the length of $\rho_n(\tau,\N)$ we use the following lemma.
\begin{lemma}
\label{lemmajogeodesique}
There exist positive constants $a_0$, $C_1$ and $C_2$ depending only on $\nu(\{0\})$ and $d$, and a random variable $Z_n$ such that if
$\nu(\{0\})<\bar{p}_c(d)$, then for every $n$ and  every $m$,
$$\PP(|\rho_n(\tau,\N)|\geq m)\leq e^{-C_1m}+\PP(Z_n> a_0m)\;,$$
and 
$$\EE(Z_n)\leq C_2n\;.$$
\end{lemma}
\begin{proof}
For any $a>0$, and $m\in\NN$,
\begin{eqnarray}
\nonumber\PP(|\rho_n(\tau,\N)|\geq m)&\leq&\PP(\exists \mbox{ s.a. }\gamma\mbox{ s.t. }|\gamma|\geq m\mbox{ and
}\sum_{e\in\gamma}\tau(e)\leq am)\\
\nonumber&+&\PP(T_n(\tau,\N)> am)\;.
\end{eqnarray}
Fix $\epsilon>0$ such that $\PP(\tau_e>\epsilon)<\bar{p}_c$ and let $X_e=\II\{\tau_e>\epsilon\}$. Then $\epsilon X_e\leq \tau_e$, and consequently
$$\PP(\exists \mbox{ s.a. }\gamma\mbox{ s.t. }|\gamma|\geq m\mbox{ and
}\sum_{e\in\gamma}\tau(e)\leq am)\leq\PP(\exists \mbox{ s.a. }\gamma\mbox{ s.t. }|\gamma|\geq m\mbox{ and
}\sum_{e\in\gamma} X_e\leq a'm)\,,$$ 
where $a'=a/\epsilon$. Choosing $a_0=\eps a_1$ with $a_1$ as in Theorem \ref{theo:lden}, this implies that  
$$\PP(\exists \mbox{ s.a. }\gamma\mbox{ s.t. }|\gamma|\geq m\mbox{ and
}\sum_{e\in\gamma}\tau(e)\leq a_0m)\leq  e^{-a_2m}\,.$$
This is the analogue of Proposition~5.8 in \cite{Kesten84}. Finally, in Corollary 4 of \cite{Pimentel10}, it is shown that there is a particular path $\gamma_n$
from $0$ to $\vec{n}$, that is independent of $\tau$ and whose expected number of edges is of order $n$. This path is constructed by walking through neighbour cells that intersect line segment $[0,\vec{n}]$ (connecting $0$ to $\vec{n}$). Hence, the size of $\gamma_n$ is at most the number of cells intersecting $[0,\vec{n}]$, which turns to be of order $n$. Thus, denoting
$Z_n=\sum_{e\in\gamma(0,\vec{n})}\tau(e)$ the passage time along this path,
we have:
\begin{equation}
\label{eq:majoTnZn}
\PP(T_n(\tau,\N)> a_0m)\leq \PP(Z_n> a_0m)\;.
\end{equation}
Using only the fact that the edge-times possess a finite moment
of order 1, there is a constant $C_2$ such that:
$$\EE(Z_n)=\EE(\tau_e)\EE(|\gamma_n|)\leq C_2n\;.$$  
\end{proof}

\vskip 2mm\noindent {\it Proof of Theorem \ref{thmmajovariance} :} Let $F$ be a function from $\Omega\times\RR_+^{\ER}$ to $\RR$ that is measurable with respect
to the completion of $\mathcal{F}\otimes \mathcal{B}(\RR_+)^{\otimes\ER}$. If for $\pi$-almost  all $\omega$, the functions $\tau\mapsto F(\omega,\tau)$ are square integrable with respect to $\nu^{\otimes \ER}$, we define the following random variables:
\begin{eqnarray*}
\EE_{\nu}(F)&:&\left\lbrace\begin{array}{rcl}\Omega&\rightarrow &\RR\\
N&\mapsto&\int F(\N,\tau)\;d\nu^{\otimes \ER}(\tau)\end{array}\right.\\
\Var_{\nu}(F)&=&\EE_{\nu}(F^2)-\EE_{\nu}(F)^2\;,
\end{eqnarray*}
and for any function $g$ in $L^2(\Omega,\pi)$,
$$\EE_{\pi}(g):=\int{g(\N)}\;d\pi(\N)\;,$$
$$\Var_{\pi}(g):=\EE_{\pi}(g^2)-\EE_{\pi}(g)^2\;.$$
We shall use the following decomposition of the variance:
$$\Var(F)=\EE_{\pi}(\Var_{\nu}(F))+\Var_\pi(\EE_{\nu}(F))\:.$$
To show that  $\EE_{\pi}(\Var_{\nu}(T_n))=O(n)$ is now standard. Indeed,
$\pi$-almost surely (see (2.17) and (2.24) in \cite{Kesten93}):
$$\Var_{\nu}(T_n)\leq 2 M_2(\nu)\EE_\nu(\rho_n)\;.$$
Thus,
$$\EE_{\pi}(\Var_{\nu}(T_n))\leq 2 M_2(\nu)\EE (\rho_n)=O(n)\;,$$
according to Lemma \ref{lemmajogeodesique}. The harder part to bound
is $\Var_\pi(\EE_{\nu}(T_n))$. Let us denote by $F_n$ the random variable
$\EE_{\nu}(T_n)$. We want to apply Proposition
\ref{prop:PoincarePoisson} to $F_n$. First note that $F_n$ belongs to
$L^2(\Omega,\pi)$: this follows from (\ref{eq:majoTnZn}). We claim that:
\begin{equation}
\label{eq:claim1V}
\forall v\in \N,\;(F_n(\N)-F_n(\N\setminus\{v\}))_-^2\leq
4M_2(\nu)d_\N(v)^2\EE_{\nu}(\II_{v\in\rho_n(\tau,\N)})\;,
\end{equation}
and
\begin{equation}
\label{eq:claim2V}
\forall x\not \in \N,\;(F_n(\N)-F_n(\N\cup\{x\}))_-^2\leq 4M_1(\nu)^2\EE_{\nu}(\II_{x\in\Gamma(\N,\rho_n(\tau,\N))})\;.
\end{equation}
Indeed, suppose first that $v$ belongs to $\N$. If
$v\not\in\rho_n(\tau,\N)$, then $\rho_n(\tau,\N)$ is still included in
$\mathcal{D}(\N\setminus\{v\})$. This implies that
$T_n(\tau,\N\setminus\{v\})\leq T_n(\tau,\N)$. Suppose that on the contrary,
$v\in\rho_n(\tau,\N)$. Let $S_1(\N,v)$ be a spanning tree of the (connected, since $d\geq 2$) subgraph induced by all the neighbors of $v$. Define a set of edges $S_2(\N,v)$ containing all the
edges of $\rho_n(\tau,\N)$ that are still in
$\mathcal{D}(\N\setminus\{v\})$, and all the edges of $S_1(\N,v)$. Then,
$S_2(\N,v)$ is a set of edges in $\mathcal{D}(\N\setminus\{v\})$ which
contains a path from $0$ to $\vec{n}$. From these considerations, we
deduce that for $v$ in $\N$:
\begin{eqnarray*}
(F_n(\N)-F_n(\N\setminus\{v\}))_-&\leq&\EE_{\nu}[(T_n(\tau,\N)-T_n(\tau,\N\setminus\{v\}))_-]\;,\\
&\leq&\EE_{\nu}\left[\left(\sum_{e\in S_2(\N,v)}\tau(e)-\sum_{e\in\rho_n(\tau,\N)}\tau(e)\right)\II_{v\in\rho_n(\tau,\N)}\right]\;,\\
&\leq&\EE_{\nu}[\sum_{e\in S_1(\N,v)}\tau(e)\II_{v\in\rho_n(\tau,\N)}]\;,\\
&\leq&\sqrt{\EE_{\nu}\left[\left(\sum_{e\in S_1(\N,v)}\tau(e)\right)^2\right]}\sqrt{\EE_\nu(\II_{v\in\rho_n(\tau,\N)})}\;,
\end{eqnarray*}
where we used Cauchy-Schwarz inequality. Using the independence of the variables $\tau(e)$ and the fact that the number of edges in $S_1(\N,v)$ is $d_\N(v)-1$, we obtain claim (\ref{eq:claim1V}). To see that claim (\ref{eq:claim2V}) is true,
suppose that $x$ does not belong to $\N$. If $x$ is not in 
$\Gamma(\N,\rho_n(\tau,\N))$, obviously $\rho_n(\tau,\N)$ is still included in
$\mathcal{D}(\N\cup\{x\})$. Thus, $T_n(\tau,\N\cup\{x\})\leq
T_n(\tau,\N)$. On the contrary, if $x$ is in
$\Gamma(\N,\rho_n(\tau,\N))$,  there are two special neighbors of $x$,
$v_{in}(x)$ and $v_{out}(x)$ such that $\rho_n(\tau,\N)$ still connects
$v_{in}(x)$ to $0$ and $v_{out}(x)$ to $\vec{n}$ in
$\mathcal{D}(\N\cup\{x\})$. Let $S_3(\N,x)$ be the set of edges
containing all the
edges of $\rho_n(\tau,\N)$ that are still in
$\mathcal{D}(\N\cup\{x\})$, plus the two edges $(v_{in}(x),x)$ and
$(x,v_{out}(x))$. $S_3(\N,x)$ contains a path in
$\mathcal{D}(\N\cup\{x\})$ from $0$ to $\vec{n}$. Notice that $v_{in}(x)$
and $v_{out}(x)$ depend on $\N$ and
$(\tau(e))_{e\in\mathcal{D}(\N)}$, and that conditionnally on $\N$, $\tau(v_{in}(x),x)$ and
$\tau(x,v_{out}(x))$ are independent from $(\tau(e))_{e\in\mathcal{D}(\N)}$. From this, we deduce that for $x$ not in
$\N$,
\begin{eqnarray*}
(F_n(\N)&-&F_n(\N\cup\{x\}))_-\\&\leq&\EE_{\nu}[(T_n(\tau,\N)-T_n(\tau,\N\cup\{x\}))_-]\;,\\
&\leq&\EE_{\nu}\Big[\big\{\sum_{e\in S_3(\N,x)}\tau(e)-\sum_{e\in\rho_n(\tau,\N)}\tau(e)\big\}\II_{x\in\Gamma(\N,\rho_n(\tau,\N))}\Big]\;,\\
&\leq&\EE_{\nu}[\{\tau(v_{in}(x),x)+\tau(x,v_{out}(x))\}\II_{x\in\Gamma(\N,\rho_n(\tau,\N))}]\;,\\
&\leq&\EE_{\nu}[\EE_{\nu}\{\tau(v_{in}(x),x)+\tau(x,v_{out}(x))|(\tau(e))_{e\in\mathcal{D}(\N)}\}\II_{x\in\Gamma(\N,\rho_n(\tau,\N))}]\;,\\
&=&2M_1(\nu)\EE_{\nu}(\II_{x\in\Gamma(\N,\rho_n(\tau,\N))})\;,
\end{eqnarray*}
which proves claim (\ref{eq:claim2V}) via Jensen's inequality. Now, these two claims together
with Proposition \ref{prop:PoincarePoisson} applied to $F_n$ give:
\begin{equation}
\label{eq:poincareFn}
\Var(F_n)\leq
4M_2(\nu)\EE(\sum_{v\in\rho_n(\tau,\N)}d_\N(v)^2)+4M_1(\nu)^2\EE(\lambda(\Gamma(\N,\rho_n(\tau,\N)))\;.
\end{equation}
We shall conclude using Lemma \ref{lemmajogeodesique}, Proposition
\ref{prop:Gamma} and Theorem \ref{majotheo}.  Let $z_1$ be as in Theorem
\ref{majotheo}, and $a$ as in Lemma \ref{lemmajogeodesique}. Notice that $f:x\mapsto x^2$ is
$c$-nice for any $c>0$.
\begin{eqnarray*}
\EE\lbrack \sum_{v\in\rho_n(\tau,\N)}f(d_\N(v))\rbrack &=&n\int_0^\infty
\PP(\sum_{v\in\rho_n(\tau,\N)}f(d_\N(v))>nz)\;dz\;,\\
&\leq&nz_1+\int_{z_1}^\infty
\PP(\sum_{v\in\rho_n(\tau,\N)}f(d_\N(v))>nz)\;dz\;,
\end{eqnarray*}
From Lemma \ref{lemmajogeodesique} and
Theorem \ref{majotheo}, we have:
\begin{eqnarray*}
\PP(\sum_{v\in\rho_n}f(d_\N(v))>nz)&\leq&
\PP(|\rho_n|>nz/z_1)+\sum_{k=0}^{nz/z_1}\PP(F_k>nz)\;,\\
&\leq&  e^{-C_1(nz/z_1)}+\PP(Z_n> a_0nz/z_1)+\frac{nz}{z_1}e^{-c_2(c_3nz)^{1/4}}\;,
\end{eqnarray*}
and thus, using that $\EE(Z_n)=O(n)$, we get:
$$\EE\lbrack \sum_{v\in\rho_n(\tau,\N)}f(d_\N(v))\rbrack =O(n)\;.$$
The proof that $\EE(\lambda(\Gamma(\N,\rho_n(\tau,\N)))=O(n)$ is
completely similar, using Proposition
\ref{prop:Gamma} and Lemma \ref{lemmajogeodesique}. Theorem
\ref{thmmajovariance} now follows from (\ref{eq:poincareFn}).

{\hfill $\square$ \vskip 2mm \noindent} 

\subsection{Stabbing number}
\label{subsec:stabbing}

The stabbing number $\mathrm{st}_n(\mathcal{D}(\N))$ of $\mathcal{D}(\N)\cap[0,n]^d$ is defined in
\cite{AddarioberrySarkar05} as the maximum number of Delaunay cells that intersect a
single line in $\mathcal{D}(\N)\cap [0,n]^d$. In
\cite{AddarioberrySarkar05}, the following deviation result 
for the stabbing number of $\mathcal{D}(N)\cap[0;n]^d$ is announced in
Lemma 3, and credited 
to Addario-Berry, Broutin and Devroye. In fact, the precise reference is
unavailable, and it seems, according to \cite{Broutin10private}, that the proof worked only in dimension 2.
\begin{lemma}[Addario-Berry, Broutin and Devroye]
\label{lem:stabbingABD}
Fix $d\geq 1$. Then, there are constants $\kappa=\kappa(d)$, $K=K(d)$
such that:
$$\EE(\mathrm{st}_n(\mathcal{D}(\N)))\leq \kappa n\;,$$
and, for any $\alpha>0$,
$$\PP(\mathrm{st}_n(\mathcal{D}(\N))>(\kappa +\alpha)n)\leq e^{-\alpha
  n/(K\log n)}\;.$$
\end{lemma}
The importance of Lemma \ref{lem:stabbingABD} is due to the fact that it is the
essential tool to prove that simple random walk on $\mathcal{D}(\N)$ is
recurrent in $\RR^2$ and transient in $\RR^d$ for $d\geq 3$. Here, we
show that our method allows to improve Lemma \ref{lem:stabbingABD} as follows.
\begin{lemma}
\label{lem:stabbingnew}
Fix $d\geq 1$. Then, there are constants $\kappa=\kappa(d)$, $K=K(d)$
such that:
$$\EE(\mathrm{st}_n(\mathcal{D}(\N)))\leq \kappa n\;,$$
and, for any $\alpha>0$,
$$\PP(\mathrm{st}_n(\mathcal{D}(\N))>(\kappa +\alpha)n)\leq e^{-\alpha
  n}\;.$$
\end{lemma}
\begin{proof}
Since the proof is very close to the proof of
Proposition \ref{prop:Gamma}, we only sketch the main differences. We divide the boundary of $[0,n]^d$ into $2^d2^{d-1}n^{d-1}$ $(d-1)$-dimensional cubes of
$(d-1)$-volume 1, and call this collection $S_n$. For each pair of
cubes $(s_1,s_2)$ in
$S_n$, we define $\mathrm{st}(s_1,s_2)$ as the maximum number of Delaunay
cells that intersect a single line-segment going from a point in $s_1$ to a
point in $s_2$. Let $V(s_1,s_2)$ be the union of those
line-segments. Notice that $\mathrm{st}(s_1,s_2)$ is bounded from
above by the number of points which belong to the union of the Delaunay cells 
intersecting $V(s_1,s_2)$. Lemma \ref{lemgoodbox} allows us to control
those points. We obtain the following analogue of Lemma \ref{lemfrontiere}.
\begin{lemma}
\label{lemfrontierestabbing} 
\begin{eqnarray*}
\PP\lbrack \mathrm{st}(s_1,s_2)>
t\rbrack&\leq&
\sum_{i=1}^{3^d}\PP\left\lbrack\sum_{\C\in\Cl^{3r,i}(\A^{3r,i}(V(s_1,s_2)))}|\C|_\N>\frac{t}{3^d.2} \right\rbrack\\
&&+\sum_{i=1}^{3^d}\PP\left\lbrack\sup_{\phi\in \SA_n}|\A^{3r,i}(V(s_1,s_2))|_\N>
\frac{t}{3^d.2}\right\rbrack\;.
\end{eqnarray*}
\end{lemma}
Note that the cardinals of $\A^{3r,i}(V(s_1,s_2)),\;i=1,\ldots,3^d$
are of order $O(n)$, and the possible choices for the pair $(s_1,s_2)$
is of order $O(n^{2(d-1)})$. Now, when $r$ is chosen large enough (see section \ref{subsec:animals}),  Lemma \ref{lem:stabbingnew} follows from Lemmas
\ref{lemmajoanimal-1}, \ref{lemcorrelneg} and \ref{lemmajotaillecluster}.

\end{proof}

\section{Appendix}
 
\subsection{Proof of Lemma \ref{lemclusterbox}}
\label{subsec:animals}
It suffices to prove Lemma \ref{lemclusterbox} for $i=1$, so we shall
omit $i$ as a subscript or superscript. Remark first that, as already noted
 at the beginning of section \ref{comp},
$p_{r}:=\sup_i\sup_z\PP(X_\z^{r,i}=1)$ tends to $0$ as $r$ tends to
infinity. Let us choose $r_0$ in such a way that for any $r\geq r_0$,
$$p_r<p_c(\GG_d)\;,$$ 
where $p_c(\GG_d)$ is the
critical probability for site percolation on $\GG_d$. When $r\geq
r_0$, we are in the so-called ``subcritical'' phase for percolation of
clusters of $(\N,6)$-bad boxes. From now on, we fix $r$ to satisfy
$r\geq r_0$ and we shall omit $r$ as a subscript or superscript, to
shorten the notations. It is well
known that in this ``subcritical'' phase, the size of the (bad)-cluster containing a
given vertex decays exponentially (see \cite{Grimmett}, Theorem (6.75)). For instance, there is a positive
constant $c$, which depends only on $r$, such that,
\begin{equation}
\label{eqclustersubcritic} \PP(|\Cl(0)|\geq x)\leq e^{-cx}\;.
\end{equation}

Remark that there is probably an exponential
number of animals of size less than $m$, and therefore, if
$\sum_{\C\in\Cl(\A)}g(|\Adh(\C)|)$ had an
exponential moment, where $g(x)=xf(x)$, it would
be easy to bound the first summand in the right-hand term of the
inequality above. When $f(x)$ is larger than $x$,
$\sum_{\C\in\Cl(\A)}g(|\Adh(\C)|)$
surely does not have exponential moments. Therefore, we have to
refine the standard argument. This refinement is a chaining
technique essentially due to \cite{Coxetal93} and we rely heavily on
that paper.

We shall prove
(\ref{eqmajoexpocluster}), the proof of (\ref{eqmajoexpobox}) being similar
and easier. We begin by stating and proving two lemmas on site percolation.

\begin{lemma}
\label{lemmajotaillecluster} 
For all $r>r_0$ there is a constant $c_9$ (depending only on $r$) such that:
$$
\PP(|\Adh(\Cl(\0))|_\N>s)\leq\PP_{p_r}(|\Adh(\Cl(\0))|_\N>s)\leq 2e^{-c_9s}\;.
$$
where by $\PP_{p_r}$, we mean that every site is open with the same
probability 
$$p_r=\sup_i\sup_z\PP(X_\z^{r,i}=1)<p_c(\GG_d)\,.$$
\end{lemma}

\begin{proof}
First, we shall
condition on the value of $X=(X_z)$. Let $u$ be a positive real
number.
$$
\PP(|\Adh(\Cl(\0))|_\N>s|X)\leq
e^{-us}\EE(e^{u|\Adh(\Cl(\0))|_\N}|X))\;.
$$
Define:
$$\partial^{\infty} \Cl(\0)=\Adh(\Cl(\0))\setminus\Im(\Cl(\0))\;.$$
Remark that, $X$ being fixed,
$|\Adh(\Cl(\0))|_\N=|\Cl(\0)|_\N+|\partial^{\infty} \Cl(\0)|_\N$ and the two
summands are independent. Through the comparison to Lebesgue's measure
\ref{eq:comparaison}, the term $|\Cl(\0)|_\N$ is stochastically dominated by a sum of
$|\Cl(\0)|$ independent random variables, each of which is a sum of
independent random variables with a Poisson distribution of parameter
$\beta=\beta(r)$, conditioned on the fact that one of them at least
must be zero. The term $|\partial^{\infty} \Cl(\0)|_\N$ is
stochastically dominated by a sum of at most $c_d|\Cl(\0)|$ random variables ($c_d$ only depending on $d$), each of which is a sum of independent random variables with Poisson distribution, conditioned on the fact that all
of them are greater than one. Obviously, the sum of independent
random variables with a Poisson distribution, conditioned on the
fact that one of them at least must be zero is stochastically
smaller than the sum of the same number of independent random
variables with Poisson distribution, conditioned on the fact that
all of them are greater than one. If $Z$ is a Poisson random
variable with parameter $\beta$, one has:
$$
\EE(e^{uZ}|Z\geq 1)=\frac{e^{\beta(e^u-1)}-e^{-\beta}}{1-e^{-\beta}}\;.
$$
Therefore,
$$
\EE(e^{u|\Adh(\Cl(\0))|_\N}|X)\leq
\left(\frac{e^{\beta(e^u-1)}}{1-e^{-\beta}}\right)^{7|\Cl(\0)|}\;,
$$
and thus:
$$
\PP(|\Adh(\Cl(\0))|_\N>s|X)\leq
e^{-us}\left(\frac{e^{\beta(e^u-1)}}{1-e^{-\beta}}\right)^{7|\Cl(\0)|}\;.
$$
Let us denote $\phi(t)=t\log t-t+1$,
$K=\frac{1}{(1-e^{-\beta})^{\frac{1}{\beta}}}$ and
$K'=\phi^{-1}(2\log K)$. Remark that $\phi$ is a bijection from
$[1,+\infty)$ onto $[0,+\infty)$. Minimizing in $u$ gives:
\begin{eqnarray*}
\PP(|\Adh(\Cl(\0))|_\N>s|X)&\leq &\left\lbrack
  Ke^{-\phi\left(\frac{s}{7\beta|\Cl(\0)|}\right)}\right\rbrack^{7\beta|\Cl(\0)|}\II_{\frac{s}{ K'}\geq7\beta|\Cl(\0)|
 }+\II_{\frac{s}{K'}\leq7\beta|\Cl(\0)|
  }\;,\\
&\leq &
  e^{-\frac{7}{2}\beta|\Cl(\0)|\phi\left(\frac{s}{7\beta|\Cl(\0)|}\right)}\II_{\frac{s}{K'}\geq
  7\beta|\Cl(\0)| }+\II_{\frac{s}{K'}\leq 7\beta|\Cl(\0)|
  }\;.
\end{eqnarray*}
Remark that:
$$
\forall x>1,\;\phi(x)\geq \min\{x-1,\frac{(x-1)^2}{e^2-1}\}\;.
$$
Therefore, if $\frac{s}{7\beta|\Cl(\0)|}\geq
  K'>1$, denoting:
$\bar{K}=\frac{1}{2}\min\{1-\frac{1}{K'},\frac{(K'-1)^2}{2K'(e^2-1)}\}$,
\begin{equation}
\label{eqminophi}
\frac{7}{2}\beta|\Cl(\0)|\phi\left(\frac{s}{7\beta|\Cl(\0)|}\right)\geq
\bar{K}s\;.
\end{equation}
This leads to:
\begin{equation}
\label{majoprobaclustercond} \PP(|\Adh(\Cl(\0))|_\N>s|X)\leq
e^{-\bar{K}s}+\II_{\frac{s}{K'}\leq 7\beta|\Cl(\0)|}\;.
\end{equation}

Therefore, denoting $c_9=\min \{\bar{K},\frac{c}{7\beta K'}\}$, we
get the desired result (notice that $\beta$ and $K'$ depend on $r$).
\end{proof}

\begin{lemma}
\label{lemcorrelneg}
Let $\Lambda$ be a finite subset of $\ZZ^d$. If $f$ is an increasing function from $\NN$ to $[1,+\infty[$,
$$\EE(\Pi_{\C\in\Cl(\Lambda)}f(|\C|))\leq
\EE(f(|\Cl(\0)|))^{|\Lambda|}\;.$$
\end{lemma}

\begin{proof}
See Lemma 7 in \cite{Pimentel10}. 
\end{proof}

We may now complete the proof of Lemma \ref{lemclusterbox}. we fix $m$
and $n$ two integers such that $n\geq m$. First, we get rid of the
variables $|\Adh(\C)|_N$ which are greater than $q(n)$.
Remark that, denoting $\Lambda_m^d=[-m,m]^d\cap \GG_d$,
$$
\max_{\A\in\Phi_m}\sum_{\C\in\Cl(r\A)}g(|\Adh(\C)|_\N)\leq
\sum_{\z\in \Lambda_m^d}g(|\Adh(\Cl(\z))|_\N)\;.
$$
Then,
\begin{eqnarray*}
\PP(\exists\; \z\in \Lambda_m^d\mbox{ s.t.
}g(|\Adh(\Cl(\z))|_\N)>q(n))&\leq &
(2m+1)^d\PP_{p_r}(g(|\Adh(\Cl(\0))|_\N)>\gamma(m))\;,\\
&\leq &(2m+1)^d\PP_{p_r}(|\Adh(\Cl(\0))|_\N>l(q(n)))\;.\\
\end{eqnarray*}

Denote by $c_9$ the constant
of Lemma \ref{lemmajotaillecluster}. For every $m$:
\begin{equation}
\label{eqmajograndesvariables} \PP(\exists \; z\in \Lambda_m^d\mbox{
s.t. }g(|\Adh(\Cl(\z))|_\N)>q(n))\leq
(m+1)^de^{-c_9l(q(n))}\;.
\end{equation}
Now, we want to bound from above the following probability:
$$
\PP(\max_{\A\in\Phi_m}\sum_{\Cl(\0)\in\Cl(\0)(r\A)}g(|\Adh(\Cl(\0))|_\N)\II_{\forall
  \z\in \Lambda_m^d,\; g(|\Adh(\Cl(\z))|_\N)\leq q(n)})\;.
$$
First, Lemma 1 in \cite{Coxetal93} (Lemma
\ref{lemrecouvrementanimal} below) remains true in our setting,
without any modification (just remark that the condition $l\leq n$
in their lemma is in fact not needed). Then, we define the following box in $\GG_d$, centered at $\x\in\GG_d$:
$$
\Lambda(\x,l)=\{(x_1+k_1,\ldots,x_d+k_d)\in\GG_d\mbox{ s.t. }(k_1,\ldots,k_d)\in
[-l,l]^d\}\;.
$$
\begin{lemma}
\label{lemrecouvrementanimal} Let $\A$ be a lattice animal of
$\GG_d$ containing $\0$, of size $|\A|=m$ and let $1\leq l$. Then,
there exists a sequence $\x_0=0,\x_1,\ldots ,\x_h\in\GG_d$ of
$h+1\leq 1+(2m-2)/l$ points such that
$$
\A\subset\bigcup_{i=0}^h\Lambda(l\x_i,2l),\;
$$
and
$$
|\x_{i+1}-\x_i|_{\infty}\leq 1,\hspace{1cm}0\leq
i\leq h-1\;.
$$
\end{lemma}

Continuing to follow \cite{Coxetal93} , we shall use Lemma
\ref{lemrecouvrementanimal} at different ``scales'' $k$, covering a
lattice animal by $1+(2m-2)/l_k$ boxes of length $4l_k+1$. We shall
choose $l_k$ later. For any animal $\A$ and $0<L,R<\infty$ define:
$$
S(L,R;\A)=\sum_{\C\in\Cl(r\xi)}g(|\Adh(\C)|_\N)\II_{L\leq
g(|\Adh(\C)|_\N)<R}\;.
$$
Suppose that $c_0$ and $(t(n,k))_{k\leq \log_{2}q(n)}$ are positive
real numbers such that:
\begin{equation}
\label{eqcondmajosumtmk} \sum_{k\leq \log_{2}q(n)}2^kt(n,k)\leq
c_0n\;.
\end{equation}
We shall choose these numbers later. Let $c'$ be a positive real
number to be fixed later also, and define $a=1+c'c_0$,
\begin{eqnarray*}
&&\PP(\exists\; \A\in \Phi_m\mbox{ with }S(0,q(n);\A)>an)\\
&\leq& \sum_{\substack{k\geq 0\\ 2^k\leq q(n)}}\PP(\exists\; \A\in
\Phi_m\mbox{
    with }S(2^k,2^{k+1};\A)>c't(n,k)2^{k})\;.
\end{eqnarray*}
Now, fix $k$ for the time being, let $\A$ be an animal of
size at most $m$, containing $\0$, let $\x_0=\0,\x_1,\ldots,\x_h$ be
as in Lemma \ref{lemrecouvrementanimal} and define
$\Lambda_{m,k}:=\bigcup_{i\leq
    h}\Lambda (l_k \x_i,2l_k)$. Clearly,
\begin{eqnarray*}
&&S(2^k,2^{k+1};\A)\\
&&\leq 2^{k+1}\left(\mbox{number of }\C\in
\Cl\left(\Lambda_{n,k}\right)\mbox{ with
}g(|\Adh(\C)|_\N)\geq
    2^k\right)\;.
\end{eqnarray*}

\begin{lemma}\label{lemmajoprobasommeclusters}
There exists $c_{10}\in(0,\infty)$ such that for all $t\geq
3|\Lambda|e^{-c_9s}$
$$
\PP\left(\sum_{\C\in\Cl(\Lambda)}\II_{g(|\Adh(\C)|_\N)>2^k}>t\right)\leq
e^{-c_{10}t}\;.
$$
\end{lemma}

\begin{proof}
Let $L$ be a positive real number. Notice that, conditionally on $X$,
$(|\Adh(\C)|_\N)_{\C\in\Cl(\Lambda)}$ are independent. Thus,
\begin{equation}
\label{eq:dev}
\PP\left(\sum_{\C\in\Cl(\Lambda)}\II_{|\Adh(\C)|_\N>s}>t\right)\leq e^{-L
  t}\EE\left(\Pi_{\C\in\Cl(\Lambda)}\EE(e^{\lambda\II_{|\Adh(\C)|_\N>s}}|X)\right)\;.
\end{equation}
But,
\begin{eqnarray*}\EE(e^{L\II_{|\Adh(\C)|_\N>s}}|X)&=&e^\lambda\PP(|\Adh(\C)|_\N>s|X)+(1-\PP(|\Adh(\C)|_\N>s|X))\;,\\
&=&1+\PP(|\Adh(\C)|_\N>s|X)(e^L-1)\;.
\end{eqnarray*}
From equation (\ref{majoprobaclustercond}), we know that:
$$
\PP(|\Adh(\C)|_\N>s|X)\leq
e^{-\bar{K}s}+\II_{\frac{s}{K'}\leq 7\beta|\C|}\;.
$$
Remark that the right-hand side of this inequality is an increasing
function of $|\C|$. Using Lemma \ref{lemcorrelneg} in equation (\ref{eq:dev}), we deduce:
\begin{eqnarray*}
&&\PP\left(\sum_{\C\in\Cl(\Lambda)}\II_{|\Adh(\C)|_\N>s}>t\right)\\
&\leq& e^{-L
  t}\left(1+(e^L-1)\left\lbrack e^{-\bar{K}s}+\PP_{p_r}(|\Cl(\0)|>
   \frac{s}{7\beta K'})\right\rbrack\right)^{|\Lambda|}\;.\\
\end{eqnarray*}
Using inequality (\ref{eqclustersubcritic}), and recalling that
$c_9=\inf\{\bar{K},\frac{c}{7\beta K'}\}$,
\begin{eqnarray*}
\PP\left(\sum_{\C\in\Cl(\Lambda)}\II_{|\Adh(\C)|_\N>s}>t\right)&\leq&
e^{-L t}\left(1+2(e^{\lambda}-1)e^{-c_9s}\right)^{|\Lambda|}\;,\\
&\leq&
e^{-L t}e^{2|\Lambda|(e^{L}-1)e^{-c_9s}}\;,\\
\end{eqnarray*}
Now, if $t\geq 2|\Lambda|e^{-c_9s}$, minimizing over
  $L$ gives:
$$
\PP\left(\sum_{\C\in\Cl(\Lambda)}\II_{g(|\Adh(\C)|_\N)>2^k}>t\right)\leq
e^{t-t\log\frac{t}{2|\Lambda|e^{-c_9s}}-2|\Lambda|e^{-c_9s}}=e^{-2|\Lambda|e^{-c_9s}\phi\left(\frac{t}{2|\Lambda|e^{-c_9s}}\right)}\;.
$$
Suppose now that $t\geq 3|\Lambda|e^{-c_9s}$. Using the same
argument which led to (\ref{eqminophi}), we get that there exists
$c_{10}>0$, depending only on $r$, such that:
$$
\forall \;t\geq
3|\Lambda|e^{-c_9s},\;\PP\left(\sum_{\C\in\Cl(\Lambda)}\II_{g(|\Adh(\C)|_\N)>2^k}>t\right)\leq
e^{-c_{10}t}\;.
$$
\end{proof}

Remark that:
$$
|\Lambda_{m,k}|\leq \frac{2m}{l_k}(4l_k+1)^2\leq 50ml_k\;.
$$
Suppose now that
\begin{equation}
\label{eqcondtmkbis} t(n,k)\geq 150nl_ke^{-c_9l(2^k)}\;.
\end{equation}
The number of choices for $\0=\x_0,\x_1,\ldots, \x_h$ in Lemma
\ref{lemrecouvrementanimal} is at most $9^h$. Therefore,
$$
\PP(\exists \;\A\in \Phi_m\mbox{ with
}S(2^k,2^{k+1};\A)>ct(n,k)2^{k+1})\leq
9^{\frac{2m}{l_k}}e^{-c_{10}c't(n,k)}\;.
$$
In view of the last inequality and (\ref{eqmajograndesvariables}), we would be happy if we had
$$
t(n,k)\geq \bar{c}l(q(n))
$$
and
$$
c_{10}c't(n,k)\geq (2\log 9).\frac{2m}{l_k}\;.
$$
Of course, we still have to check conditions
(\ref{eqcondmajosumtmk}) and (\ref{eqcondtmkbis}), and we can choose
$c'$ as large as we need. A natural way to proceed is first to
choose $l_k$ in such a way that the right-hand side in
(\ref{eqcondtmkbis}) is proportional to $\frac{m}{l_k}$, and then to
take large enough. Recall that $n\geq m$. Thus, define:
$$
l_k=\left\lceil e^{\frac{c_9}{2}l(2^k)}\right\rceil\;.
$$
Choose
$$
t(n,k)=\max\{150ml_k e^{-c_9l(2^k)},l(q(n))\}\;.
$$
Let $c'$ be such that:
$$
c'\geq \frac{4\log 9}{150c_{10}}\;.
$$
This ensures that:
$$
c'c_{10}t(n,k)\geq ( 4\log 9)nl_ke^{-c_9l(2^k)}\geq (2\log
9)\frac{2n}{l_k}\;.
$$
Therefore,
\begin{equation}
\label{eqmajoquasifinale} \PP(\exists \;\A\in \Phi_m\mbox{ with
}S(2^k,2^{k+1};\A)>c't(n,k)2^{k+1})\leq e^{-c_{10}c't(n,k)/2}\;.
\end{equation}
Condition (\ref{eqcondtmkbis}) is trivially verified from the
definition of $t(n,k)$. Now, let us check condition
(\ref{eqcondmajosumtmk}). Now, assume that f is
$\frac{4}{c_9}$-nice. Then, using the definition of $l_k$,
$$
\sum_{k\leq
  \log_{2}q(n)}2^{k+1}150nl_k e^{-c_9l(2^k)}\leq 300 n
  \sum_{k}e^{\log (2^k)-\frac{c_9}{2}l(2^k)}\;.
$$
By our assumption,
$$
\limsup_{k\to\infty} \frac{e^{\log
    (2^k)-\frac{c_9}{2}l(2^k)}}{(2^{\frac{1}{3}})^k}\leq 1\;,
 $$
and therefore,
$$
\Sigma:=\sum_{k}e^{\log (2^k)-\frac{c_9}{2}l(2^k)}<\infty\;.
$$
On the other hand,
$$
\sum_{k\leq
  \log_{2}q(n)}2^{k+1}l(q(n))\leq
  4 q(n)l(q(n))\leq 4n\;,
$$
by definition of $q(n)$. Therefore, condition
(\ref{eqcondmajosumtmk}) is checked with
$$
c_0=300\Sigma+4\;.
$$
Remark also that $g(x)\geq x$, and hence $l(y)\leq y$.
Therefore:
$$
n\leq (q(n)+1)l(q(n)+1)\leq (q(n)+1)^{2}\;,
$$
$$
q(n)\geq m^{1/2}-1\;,
$$
and thus, if
\begin{equation}
\label{eqcondsomme2} \liminf_{y\to\infty}\frac{l(y)}{\log
y}\geq\frac{(4d+1)}{c_9}\;
\end{equation}
then
$$
d\log (m+1)-\frac{c_9}{2}l(q(n))\xrightarrow[n\to\infty]
{}-\infty\;,
$$
and therefore, there exist a constant $c_{11}$ such that:
$$
\forall n\geq m,\;(m+1)^de^{-c_9l(q(n))}\leq
c_{11}e^{-\frac{c_9}{2}l(q(n))}\;.
$$
Since $q(n)\leq n$, we get in the same way that there is a constant
$c_{12}$ such that:
$$
\sum_{k\leq
  \log_{2}q(n)}e^{-c_{10}ct(n,k)/2}\leq
  \log_{2}q(n)e^{-c_{10}cl(q(n))/2}\leq
  c_{12}e^{-c_{10}cl(q(n))/4}\;,
  $$
provided that (\ref{eqcondsomme2}) holds. Therefore, we define
  $c_5=\sup\{\frac{6}{c},\frac{4d+1}{c_9}\}$, and we suppose that
\begin{equation}
\label{eqcondsommefinale} \liminf_{y\to\infty}\frac{l(y)}{\log
y}\geq c_5\;.
\end{equation}
Then, for $c_6=a=1+cc_0$, there exists a positive constant $c_7$
such that:
\begin{eqnarray*}
\PP(\sup_{\A\in\Phi_m}\sum_{\C\in\Cl(
  r\xi)}g(|\Adh(\Cl(\0))|_\N)> c_6n)&\leq& (m+1)^2e^{-\bar{c}l(q(n))}\\ &+&\sum_{k\leq
  \log_{2}q(n)}e^{-c_{10}ct(n,k)/2}\;,\\
&\leq &e^{-c_7l(q(n))}\;,
\end{eqnarray*}
where the first inequality follows from (\ref{eqmajograndesvariables}) and (\ref{eqmajoquasifinale}). This concludes the proof of Lemma \ref{lemclusterbox}.
\hfill $\square$

\subsection{Two simple geometric lemmas}
\label{subsec:lemchemborne}
\begin{lemma}
\label{lem:cheminborne}
Let $\N$ be a locally finite set of points in $\RR^d$ and $\D$ the
Delaunay triangulation based on $\N$. Let $u$ and $v$ be two distinct points in
$\N$. Then, there is a path on $\D$ going from $v$ to $u$ and totally
included in the (closed)  ball of center $u$ and radius $|u-v|$.
\end{lemma}
\begin{proof}
First we show that there is a neighbour $w$ of $v$ in $\D$ inside the
(closed) ball of diameter $[u,v]$. To see this, let us define by
$(B_\alpha)_{\alpha\in[0,1]}$ the collection of euclidean balls such
that $B_{\alpha}$ has diameter $[v,x_\alpha]$, where
$x_\alpha=v+\alpha (u-v)$. Notice that $B_\alpha\subset B_{\alpha'}$
as soon as $\alpha \leq \alpha'$. Define:
$$\alpha_0:=\min\{\alpha\in[0,1]\mbox{ s.t. }\exists w'\in
\N\setminus\{v\}\cap B_\alpha\}\;.$$
This is indeed a minimum because $\N$ is locally finite. Notice also
that the set is non-empty since it contains $\alpha=1$. Thus the
interior of $B_{\alpha_0}$ does not intersect $\N$, but $\partial
B_{\alpha_0}$ contains (at least) two points of $\N$ (including $v$). This implies that
there is a point on the sphere $\partial B_{\alpha_0}$ which is a
neighbour of $v$. So we have proved that there is a neighbour $w_1$ of $v$
inside the ball of diameter $[u,v]$. Then, as long as the neighbour
obtained is different from $u$, we may iterate this
construction to get a sequence of neighbours $w_0=v,w_1,w_2,\ldots$
such that $w_{i+1}$ belongs to the ball of diameter $[u,w_i]$. All
these balls are included in the ball of center $u$ and radius
$|u-v|$. Since $\N$ is locally finite, this construction has to stop
at some $k$, when the condition that $w_k$ is distinct from $u$ is no
longer satisfied. Then the desired path is constructed.
\end{proof}
For the next lemma, recall that a locally finite set of points $\N$ in $\RR^d$ is said to be ``in generic position'' if every subset of points of cardinal $d+1$ can be circumscribed a unique $d$-dimensional sphere and if no $d+2$ points in $\N$ are co-spherical, i.e. lie on a common sphere in $\RR^d$. It is well known that a Poisson random set with intensity comparable to the Lebesgue measure is in generic position with probability 1.
\begin{lemma}
\label{lemm:gammaNe}
Let $\N$ be a locally finite set of points in $\RR^d$ in generic position and define, for any edge $e$ of the Delaunay graph $\D(\N)$:
$$\Gamma(\N,e)=\{x\in\RR^d\mbox{ s.t. }e\not\subset
\mathcal{D}(\N\cup \{x\}\}\;.$$
Then,
$$\Gamma(\N,e)=\bigcap_{\Delta\ni e}\inter{B}(\Delta)\;,$$
where the intersection is taken  over all Delaunay cells $\Delta$
which contain $e$.
\end{lemma}
\begin{proof}
 Let us define, for any pair of vertices $\{u,v\}$ in $\N$, the following convex $(d-1)$-dimensional polytope:
$$P(\{u,v\}):=\{y\in\RR^d\mbox{ s.t. }d(y,u)=d(y,v)\leq d(y,w)\;\forall w\in\N\}\;,$$
where $d$ is the euclidean distance. $P(\{u,v\})$ is non-empty if and only if $C_u$ and $C_v$ share a $(d-1)$-dimensional face i.e when $\{u,v\}$ is an edge of $\D(\N)$, and in this case, $P(\{u,v\})$ is precisely this common face. Now, let $e=\{u,v\}$ be an edge of $\D(\N)$, and let $x$ be a point of $\RR^d$. Then, $e\not\subset\D(\N\cup\{x\})$ if and only if all the points of $P(e)$ are (strictly) closer to $x$ than to $u$ (or $v$, but this is the same since $d(y,u)=d(y,v)$ for $y$ in $P(e)$). But this is a convex condition: this merely means that $P(e)$ is included in the open half space containing $x$ and with boundary the median hyperplane of $[x,u]$. Thus $e\not\subset\D(\N\cup\{x\})$ if and only if all the extreme points of the convex $(d-1)$-dimensional polytope $P(e)$ are closer to $x$ than to $u$. But these extreme points are those points $y\in P(e)$ for which there are $d-1$ different points $w$ in $\N\setminus\{u,v\}$ such that $d(y,u)=d(y,v)=d(y,w)$ (recall that $\N$ is in generic position). These are the centers of the circumballs of the Delaunay cells. Thus $e\not\subset\D(\N\cup\{x\})$ if and only $x$ belongs to  $\bigcap_{\Delta\ni e}\inter{B}(\Delta)$.
\end{proof}

\section*{Acknowledgements}
We would like to thank Pierre Calka and Marie Th\'eret for  helpful discussions. We are also grateful to an anonymous referee for his careful reading and relevant remarks which helped us to improve this paper.


\end{document}